\documentclass[12pt]{article} 
\usepackage{lastpage}
\newcommand{\pyear}{2024}  
\newcommand{\vol}{22}  
\newcommand{\no}{3} 
\newcommand{\stpage}{1}  

\usepackage[top=2.54cm, bottom=2.54cm, left=2.54cm, right=2.54cm]{geometry} 
\usepackage{times} 
\usepackage{lmodern} 
\usepackage{hanging} 
\usepackage{graphicx} 
\usepackage[english]{babel}  
\usepackage{natbib}
\usepackage{fancyhdr} 
\fancypagestyle{frontpagefooter}{%
  \fancyhead{} 
  \fancyfoot[L]{\footnotesize 
    Corresponding Author: Xiao-Li Meng \\
    Email: meng@stat.harvard.edu}
  \fancyfoot[C]{\footnotesize }	  
  \fancyfoot[R]{\footnotesize } 
}
\newcommand{\E}{{\rm E}}
\newcommand{\Cov}{{\rm Cov}}
\newcommand{\Var}{{\rm V}}
\newcommand{\Corr}{{\rm Corr}}
\newcommand{\iidsim}{\overset{\text{iid}}{\sim}}
\usepackage{amsthm}

\newtheoremstyle{sa}
{3mm} 
{2mm} 
{} 
{} 
{\bfseries} 
{:} 
{.5em} 
{} 

\theoremstyle{sa}
\newtheorem{theorem}{Theorem}
\newtheorem{corollary}{Corollary}

\renewenvironment{proof}{{\noindent \bfseries Proof:}}{\qed}

\usepackage{amsmath,amssymb}
\usepackage{hyperref} 
\usepackage{lscape}
\usepackage{arydshln}
\usepackage{booktabs}
\usepackage{verbatim}
\usepackage{enumitem}
\usepackage{bm}
\DeclareMathAlphabet{\mathpzc}{OT1}{pzc}{m}{it}
\usepackage{wasysym}
\usepackage{esint}
\usepackage{xcolor}
\usepackage{tikz}
\usetikzlibrary{shapes,backgrounds,arrows,chains,matrix,positioning,scopes}
\usepackage{caption}
\usepackage{titlesec}

\titleformat{\section}{\bf}
{\makebox[1.27cm][l]{\thesection.}}{0in}{} 
\titlespacing*{\section}{0pt}{3mm}{0pt}

\titleformat{\subsection}{\bf}
{\makebox[1.27cm][l]{\thesubsection.}}{0in}{}  
\titlespacing*{\subsection}{0pt}{3mm}{0pt}

\titleformat{\subsubsection}{\bf}
{\makebox[1.27cm][l]{\thesubsubsection.}}{0in}{}  
\titlespacing*{\subsubsection}{0pt}{3mm}{0pt}

\usepackage{indentfirst}  
\usepackage[labelfont=bf,textfont=bf]{caption}  
\usepackage{lineno}

\begin{document}
\thispagestyle{empty}

\thispagestyle{frontpagefooter}
 
 \begin{minipage}[t]{0.65\textwidth}
 	\vspace{0pt}
 	\begin{flushleft}
 		For a Special Issue of \textit{Statistics and Applications} (\href{http://www.ssca.org.in/journal.html}{{{http://www.ssca.org.in/journal}}}) in Memory of C R Rao 
 	\end{flushleft}
 \end{minipage}

%
%

\begin{center}
{\fontsize{16}{24} \selectfont 
{\bf A Heisenberg-esque  Uncertainty Principle for \\ Simultaneous (Machine) Learning and Error Assessment?} 
 } \\[6mm]
\normalsize
{\bf Xiao-Li Meng \\ 
{\it $^1$Department of Statistics, Harvard University
}}

\medskip
Received: 28 September 2024; Revised: 31 November, 2024

\end{center}
\vspace{-6mm}

\setlength\parindent{1.27cm}

\noindent\rule{16.5cm}{0.8pt}\vspace{-3mm}

\noindent{\bf Abstract}\vspace{-1mm}

A highly cited and inspiring article by \cite{bates2024cross} demonstrates that the prediction errors estimated through cross-validation, Bootstrap or Mallow’s $C_P$ can all be independent of the actual prediction errors. This essay hypothesizes that these occurrences signify a broader, Heisenberg-like uncertainty principle for learning: optimizing learning and assessing actual errors using the same data are fundamentally at odds. Only suboptimal learning preserves untapped information for actual error assessments, and vice versa, reinforcing the `no free lunch' principle. To substantiate this intuition,  a Cramér-Rao-style lower bound is established under the squared loss, which shows that the relative regret in learning is bounded below by the square of the
correlation between any unbiased error assessor and the actual learning error. Readers are invited to explore generalizations, develop variations, or even uncover genuine `free lunches.' The connection with the Heisenberg uncertainty principle is more than metaphorical, because both share an essence of the Cramér-Rao inequality: marginal variations cannot manifest individually to arbitrary degrees when their underlying co-variation is constrained, whether the co-variation is about individual states or their generating mechanisms, as in the quantum realm. A practical takeaway of such a learning principle is that it may be prudent to reserve some information specifically for error assessment rather than pursue full optimization in learning, particularly when intentional randomness is introduced to mitigate overfitting.


 \vspace{-5mm} 

\noindent\rule{16.5cm}{0.8pt} 
\setcounter{equation}{0}


\section{A Rao-esque apology and a quantum-leap excuse }\label{sec:Apology}

Many of the advances in statistics and machine learning is about using data as efficiently and reliably as possible to achieve a host of learning objectives, such as inference, prediction, classification, etc.  Being statistically efficient typically means to optimize over some criterion that amounts to minimizing learning errors based on data at hand, whether in a brute-force fashion, such as minimizing a $\chi^2$ distance or adopting the $L^2$-loss directly on the target of learning, or through deeper principles, \textit{e.g.}, by maximizing a likelihood function or a posterior density. Since the actual learning errors themselves cannot be known without an external benchmark, we seek clever and reliable ways to assess them, whether for training machine learning algorithms, constructing confidence intervals, or checking Bayesian models. 

Naturally, we wish to be able to optimally use our data for both purposes: to most efficiently learn whatever we can learn, and to most reliably assess the errors in whatever we cannot learn.  However, since any information on the actual learning error can be used to improve the learning itself, we should be mindful  that optimizing one endeavor comes at the expense of the other. To emphasize this no-free lunch principle, this essay first revisits seemingly quaint examples and classical results to remind ourselves that this principle has been in action for as long as statistical inference exists. However,   such an issue  has not received much emphasis apparently because principled statistical methods, such as likelihood or Bayesian methods, automatically prioritize optimal learning over error assessment. 

Yet time has changed. Machine learning and other pattern-seeking methods require much intuition and judgment to tune well, when their theoretical guiding principles are not well developed or digested. Substituting---not merely supplementing---virtual trials and errors for sapient contemplation and introspection is becoming increasingly habitual, making us more vulnerable to wishful thinking, misinformed intuitions, and misguided common sense. To better prepare students and newcomers to our progressively empiricism-slanted culture of learning, this essay then recasts a classical result regarding UMVUE to the broader class of problems of unbiased learning, and establishes a mathematical inequality that captures the aforementioned Heisenberg-esque uncertainty principle for simultaneous learning and error assessment under the squared loss. 

This inequality is a low-hanging fruit in establishing a general theory for understanding the competing nature between optimal learning and actual error assessing.  Nevertheless, it can help us anticipate and better appreciate further results such as those obtained in \cite{bates2024cross}, which show that the error estimates from cross validations and other popular methods can be independent of actual learning error. The uncertainty principle tells us that this should not come as a surprise. Rather, the independence is an indication that the corresponding learning is optimal in some sense. 

Since this essay was prepared for this special issue in memory of Professor C. R. Rao,  it seems fitting to quote \cite{rao1962efficient}, a discussion article presented\footnote{As a reminder of C. R. Rao's remarkable longevity of life and professional life, this presentation took place before my parents had decided to conceive me.}  to the Royal Statistical Society in England  (RSS):

\begin{quote}
``While thanking the Royal Statistical Society for giving me an opportunity to read a paper at one of its meetings, I must apologize for choosing a subject which may appear somewhat classical. But I hope this small attempt intended to state in precise terms what can be claimed about m.l. estimates, in large samples, will at least throw some light on current controversies."
\end{quote}

\cite{rao1962efficient} was a paper on ``Efficient estimates and optimum inference procedures in large samples" (and his ``m.l." referred to maximum likelihood, not machine learning), one of a series of fundamental articles he authored during what is now considered an era of classical mathematical statistics. Therefore, initially I was somewhat surprised by Rao's apologetic sentiment---one that I ought to adopt myself for bringing up UMVUE in an era where few statistics students would recognize the acronym without Googling it. However, upon reflection, and considering his training under R. A. Fisher and the characteristically wry culture of RSS discussion at that time, I suspect Rao's apology was more of a gentle reminder to not ignore established literature or wisdom when facing new problems. I am therefore grateful to the editors of this special issue, especially Bhramar Mukherjee, for the opportunity to honor Professor C. R. Rao with one more example of the value of such a reminder:  how classical statistical results can offer insights and contextualization for modern work in data science like \cite{bates2024cross}.  

I am also deeply grateful to Bhramar for her extraordinary patience in allowing me two extra months to complete this essay, without which I would have embarrassed myself significantly more by writing about Heisenberg Uncertainty Principle (HUP) while knowing almost surely nothing even about classic mechanics\footnote{Majoring in pure math in 1980s China means that I had taken no courses outside of mathematics, with the exception of mandatory ones for regulating students' bodies or minds.}. The connection between Cram\'er-Rao inequality and HUP has long been suspected, but I was unaware of any statistical literature on the connection between the two (however, during this work, I was made aware of such results in information theory---see Section~\ref{sec:definingcov}).  

Unfortunately, I had found neither the time nor the courage to explore quantum physics. Bhramar’s invitation gave me a great excuse to delve into it, though clearly it has been a quantum leap (or dive). I am therefore deeply grateful to the physicists, philosophers, and statisticians (see acknowledgment) who generously took the time to educate and inspire me, 
introducing me to numerous articles that, no doubt, will require another quantum-leap excuse to digest fully. These include physics literature on quantum Cramér-Rao bounds and quantum Fisher information \citep[\textit{e.g.},][]{toth2013extremal,toth2022uncertainty}, as well as statistical writings on the relevance of quantum uncertainty to statistics \citep[\textit{e.g.},][]{gelman2013does}, to name just a few.

Nevertheless, to set readers’ expectations realistically, this essay offers nothing about the Heisenberg Uncertainty Principle (HUP) that isn’t already in Wikipedia. I wrote much of it as reading notes to educate myself, so, paraphrasing a most memorable chiasmus from an RSS discussion: “The parts of the paper that are true are not new, and parts that are new are not true” \citep{mccullagh1999discussion}. My hope, however, is that these notes may still be of use to those who share my curiosity (and innocence). I also hope that my attempt to extend the notion of covariance to quantum operators might encourage us to step out of our comfort zones without stepping out of our minds.
 
Intellectually, quantum indeterminacy is a captivating and challenging topic, especially for those of us who have been probability-law abiding citizens.  To my knowledge, currently only a few statisticians—most notably Richard Gill\footnote{See https://www.math.leidenuniv.nl/~gillrd/}—have studied it systematically. Therefore, even if everything ``new'' in this essay ends up merely demonstrating that humans can out-hallucinate ChatGPT, I’d still be content dedicating it to the legendary C. R. Rao. Throughout his extraordinary career, Professor Rao applied his statistical insight and mathematical skills to establish and solidify the foundations of statistics. As quantum computing looms on the horizon, some statisticians should be leading the way in building the foundations of quantum data science, as articulated in the discussion article ``When Quantum Computation Meets Data Science: Making Data Science Quantum'' by \cite{Wang2022When}, a prominent statistician exploring quantum computing’s role in data science. Thus, even if this essay inspires only one future C. R. Rao of quantum data science, it won’t take a quantum leap to believe that Professor Rao would embrace my dedication.

More broadly, I would find great professional satisfaction (and justification for my insomnia) if this essay serves as a reminder that time-honored statistical theory and wisdom have much to offer as we statisticians are increasingly called to step outside our comfort zones—from embracing machine learning to anticipating quantum computing. By learning from and contributing to other fields, especially time-tested ones such as philosophy and physics, we can enhance the intellectual impact of our discipline.

\section{A paradox of error assessment?}\label{sec:paradoz}

Let us start with an excursion to the classical statistical sanctuary most frequently adopted in statistical research and pedagogy: we have an independently and identically distributed (i.i.d) normal sample,  $X_1, \cdots, X_n \iidsim N(\mu, \sigma^2)$, and we are interested in making inference about $\mu$. It is well-known that the maximum likelihood estimator (MLE) for $\mu$ is the sample mean $\bar X_n$.  The actual error of the MLE then is  $\delta=\bar X_n - \mu$.  It is textbook knowledge that the sample mean $\bar X_n$ and the sample variance $S_n^2$ are independent under the normal model $N(\mu, \sigma^2)$.  This fact is critical for establishing perhaps the most celebrated pivotal quantity in statistics, $t=\sqrt{n}(\bar X - \mu)/S_n$, \textit{i.e.}, the $t$ statistic, because of the existence of parameter-free distribution of $t$ for any $n\ge 2$,  thanks to the aforementioned independence.

But this independence also implies a seemingly paradoxical fact that has received no mention in any textbook (that I am aware of): that $\hat\delta^2 \equiv S_n^2/n$ apparently is the worst estimate of the square of the actual error  $\delta^2 = (\bar X_n - \mu)^2$, because $\hat\delta^2$ and $\delta^2$ are independent of each other for any choice of $\theta=\{\mu, \sigma^2\}$. In what other context would a statistician (knowingly) suggest estimating an unknown with an independent quantity?

The article by \cite{bates2024cross} reminds us that this seemingly paradoxical phenomenon is far more prevalent than we may have realized. To recast their findings in a broader setting but with a scalar estimand for notational simplicity, consider the possibly heteroscedastic linear regression setting, 
\begin{equation}\label{eq:hetroregre}
Y_i=\theta X_i + \epsilon_i,  \quad \text{where} \quad \E[\epsilon_i|\mathbf X]=0, \Var(\epsilon_i|\mathbf X)=\sigma^2_i, \quad i=1, \ldots, n.
\end{equation}
and conditioning on $\mathbf{X}=\{X_1, \ldots, X_n\}$,  $\{\epsilon_1, \ldots, \epsilon_n\}$ are mutually independent.  As \cite{bates2024cross} reminds us,  when $\{\epsilon_1, \ldots, \epsilon_n\}$ are i.i.d $N(0, \sigma^2)$,  the least-squares estimator for $\theta$, $
\hat\theta_{\rm LS} = \sum_{i=1}^n Y_iX_i/{\sum_{i=1}^n X_i^2}$ 
is independent of the residual $R=\{\hat r_i=Y_i - \hat\theta X_i, i=1, \ldots, n\}$, for any given $\{\theta, \sigma^2\}$. Consequently, since the true predictive error depends on the data only through $\hat\theta_{\rm LS}$, and cross-validation error estimators are functions only of the residuals, the true and estimated errors are independent of each other. The results obviously apply to any error estimates that depend on data only through $R$, which are the case virtually for all the common estimators in practice, as demonstrated in \cite{bates2024cross}.

It is well-known \cite[\textit{e.g.},][]{casella2024statistical} that under the i.i.d normal setting, $\hat\theta_{\rm LS}$ is the MLE and indeed UMVUE (uniformly minimum variance unbiased estimator) because its variance reaches the Cram\'er-Rao bound.  Even without the normality, we know that $\hat\theta_{\rm LS}$ is BLUE (best linear unbiased estimator) and it is linearly uncorrelated with the residual $R$ under the squared loss, because it is the orthogonal projection of $Y$ onto the space expanded by $\mathbf X$ when $\sigma_i$ is invariant of $i$. 

Although rarely mentioned in textbooks, this optimality-orthogonality duality appears in essentially all inferential paradigms.  Geometrically speaking, the equivalence is due to the fact that the linear correlation between two variables is the cosine of the angle between them in the $L^2$ space, and optimal projection is the orthogonal projection. Probabilistically, the ubiquity of this duality is manifested  by the so-called ``Eve's law" \citep{blitzstein_hwang_2014},  an instance of the Pythagoras theorem in the $L^2$ space. 

That is, under any joint distribution, $p(H, G)$, as long as it generates finite second moments, $\Cov[H-\E(H|G), \E(H|G)]=0$, because $\E(H|G)$ is the orthogonal projection of $H$ to the space of $L^2$ functions that are measurable with respect to the $\sigma$-field generated by $G$. Consequently, the Pythagoras theorem is in force:
\begin{align}\label{eq:eve}
\Var(H)&=\E\left[H-\E(H)\right]^2 \nonumber
=\E[H-\E(H|G)]^2+\E[\E(H|G)- \E(H)]^2\nonumber\\
&=\E[\Var(H|G)]+\Var[\E(H|G)], 
\end{align}
which is Eve's law.  The ubiquity of the duality is due to the fact that the expectation operator in \eqref{eq:eve} can be taken with any kind of distribution: posterior (predictive) distributions for Bayesian inferences, super-population distributions as typical for likelihood inference (as in the $N(\mu, \sigma^2)$ example), or randomization distributions as in finite-population calculations \citep[as adopted in][]{meng2018statistical}.  

Nevertheless, this duality is a qualitative statement, as it does not quantify what happens for non-optimal estimation or learning.  As demonstrated below, this duality can be extended quantitatively by tethering the deficiency in learning with the relevancy in assessing the actual learning errors. This quantification crystallizes the reason for the apparent paradox, and it can help reduce wasted efforts in pursuit of the impossible. It also makes it clearer that there is no real paradox, much like how Simpson's paradox is not a paradox once its workings are revealed and understood \citep[\textit{e.g.}][]{liu2014comment,gong2021judicious}. 

The title of the next section says it all: there is no free lunch. If there is any data information left—after learning—for assessing the actual error, then we can reduce the actual error by removing the part that can be predicted by the untapped data information. This implies our learning is not optimal, and vice versa. 
Section~\ref{sec:lunch} illustrates this fact in the context of heteroscedastic regression, followed by a broad reflection in Section~\ref{sec:Jay} on its implications in the context of error assessment without external benchmarks, a statistical magic. Sections~\ref{sec:mainresult} and \ref{sec:asymresults} then establish respectively the exact and asymptotic inequalities that capture the learning uncertainty principle under the squared loss. 

To facilitate a formal comparison with the Heisenberg Uncertainty Principle (HUP) using the notion of co-variation, Section~\ref{sec:definingcov} discusses the generalization of the measure of co-variance from real-valued variables to complex-valued variables and functions. Section~\ref{sec:covoperator} then applies the generalization to the case of HUP by defining co-variances between mechanisms (\textit{e.g.}, the position and momentum \textit{operators}) rather than between the states they generate (\textit{e.g.}, the actual position and momentum states). 
With these preparations, Section~\ref{sec:covariation} compares the learning-error inequality, Cram\'er-Rao inequality, and HUP inequality, highlighting their shared essence from a statistical perspective.  

Section~\ref{sec:philosophy} reflects on various philosophical issues surrounding uncertainty principles in general, and HUP in particular, with insights from the encyclopedic essay by \cite{hilgevoord2024uncertainty}. Section~\ref{sec:marriage} briefly touches on the trade-off between quantitative and qualitative studies, prompted by a discussion in \cite{hilgevoord2024uncertainty}, and how intercultural inquires can benefit from their happy marriage. This leads to a piece of advice from Professor Rao on living a happy life, which serves as a fitting conclusion to this essay in his memory. However, to encourage students to engage with this essay to the fullest extent of their attention spans, Section~\ref{sec:prologue} provides a prologue, especially for those who may not enjoy technical appendices but wish the essay were even longer.

\section{Once again, there is no free lunch}\label{sec:lunch}

Consider the heteroscedastic setting \eqref{eq:hetroregre}, where we know that BLUE is given by the weighted LS, in the form of
\begin{equation}\label{eq:wSL}
\hat\theta_{w} = \frac{\sum_{i=1}^n w_iY_iX_i}{\sum_{i=1}^n w_iX_i^2},
\end{equation}
when the weights $w_i \propto \sigma^{-2}_i, i=1, \ldots, n$.  Now consider an arbitrarily weighted $\hat\theta_w$, and its correlation---denoted by $\rho$---with the corresponding residual $R_w=\{\hat r_{w, i}=Y_i - \hat\theta_w X_i; i=1, \ldots, n\}$. For conveying the main idea, the case of $n=2$ is sufficient.  As a special case of the general expression given in Appendix A, we have, conditioning on $\mathbf X$ (but we suppress this conditioning notation-wise unless necessary), 
\begin{equation}\label{eq:corrn2}
\rho^2(\hat\theta_w, \hat r_{w, i}) =
\frac{X_1^2X_2^2(w_1\sigma_1\sigma_2^{-1}-w_2\sigma_2\sigma_1^{-1})^2}{(w_1^2X_1^2 \sigma^2_1+w^2_2X_2^2 \sigma^2_2)(X_1^2\sigma^{-2}_1+ X_2^2\sigma^{-2}_2)}, \quad i=1, 2,
\end{equation}
which is zero if and only if $w_i \propto \sigma^{-2}_i, i=1, 2$ (as long as $X_i\not=0, i=1, 2$). That is, $\hat\theta_w$ is BLUE
(or the MLE if we assume normality) if and only if $\hat\theta_w$ is uncorrelated with $\hat r_{w,i}$. More importantly, expression \eqref{eq:corrn2} tells us exactly how the statistical efficiency of $\hat\theta_w$ is directly linked to this correlation.  

Specifically, let $\hat\theta_{\rm BLUE}$ be the optimally weighted LS estimator with weight $w_i \propto \sigma^{-2}_i, i=1, 2$, and $RR_w$ be the relative regret of an arbitrarily weighted $\hat\theta_w$ under the squared loss, that is,
\begin{align}\label{eq:regret}
RR_w = \frac{\Var(\hat\theta_w) - \Var(\hat\theta_{\rm BLUE})}{\Var(\hat\theta_w)}
= 1- \frac{(w_1X_1^2+w_2X_2^2)^2}{(w_1^2X_1^2 \sigma^2_1+w_2^2X_2^2 \sigma^2_2)(X_1^2 \sigma^{-2}_1+X_2^2 \sigma^{-2}_2)}.
\end{align}
Whereas it may not be immediate from \eqref{eq:corrn2} and \eqref{eq:regret}, one can verify directly that
\begin{equation}\label{eq:corrrr}
\rho^2(\hat\theta_w, \hat r_{w, i}) = RR_w,  \quad i=1, 2,
\end{equation}
for any choice of weights $w$ or values of $\{\sigma^2_i, i=1, 2\}$. This means that if we want to increase the magnitude of the correlation between $\hat\theta_w$ and $\hat r_{w, i}$, we must sacrifice the efficiency of $\hat\theta_w$, and vice versa. 

But why would we want to increase $|\rho(\hat\theta_w, \hat r_{w, i})|$?   Consider the case where our learning target is $c\theta$, with $c$ being a constant. For example,  we take $c=1$ when the regression coefficient $\theta$ is the target, or $c=X^*$ when the learning target is the mean of $Y$ when $X=X^*$. In such cases, the actual error is given by $\delta_w=c(\hat\theta_w -\theta)$. We can assess $\delta_w$ via $\hat\delta_w=\tilde c\hat{r}_{w,1}$ for some choice of $\tilde c$  (recall $\hat r_{w,1}+\hat r_{w,2}=0$ and hence a single residual suffices).  Because 
\begin{align}\label{eq:corr3}
\rho^2(\delta_w, \hat\delta_w) =\rho^2(c\hat\theta_w, \tilde c \hat{r}_{w,1})= \rho^2(\hat\theta_w, \hat{r}_{w,1}),
\end{align}
we see that by moving $\rho^2(\hat\theta_w, \hat{r}_{w,1})$ away from zero, we will have an assessment $\hat\delta_w$ of the actual error $\delta_w$ that has some degrees of conditional relevancy, that is, $\hat\delta_w$ is at least correlated with $\delta_w$ conditioning on the setting \eqref{eq:hetroregre}. But this gain of relevancy is achieved necessarily by increasing the relative regret (recall the relative regret for $c\hat\theta_w$ is invariant to the value of $c$), that is, by sacrificing the efficiency of $\hat\theta_w$, because 
\begin{align}\label{eq:corr4}
\rho^2(\delta_w, \hat\delta_w) = RR_w,
\end{align}
thanks to \eqref{eq:corrrr}-\eqref{eq:corr3}. 

If our learning target is to predict (a new) $Y^*$ when $X=X^*$, then the actual prediction error is $\delta_w^* = Y^* - \hat\theta_w X^*$. In such cases, the prediction risk under the squared loss is
\begin{align*}
    \E(Y^* - \hat\theta_w X^*)^2= \Var(Y^*)+ (X^*)^2\Var(\hat\theta_w).
\end{align*}
Because $\Var(Y^*)$ and $(X^*)^2$ are invariant to the weights, we obtain the relative regret for prediction $RR^*_w=\gamma RR_w$, where $RR_w$ is from \eqref{eq:regret} and the adjustment factor $\gamma$ is given by
\begin{align}\label{eq:adjust}
    \gamma=\frac{(X^*)^2\Var(\hat\theta_w)}{ \Var(Y^*)+ (X^*)^2\Var(\hat\theta_w)}.
\end{align}
Furthermore, because $\hat\delta_w=\tilde c \hat r_{w,1}$ is independent of $Y^*$,
$\Cov(\delta_w^*, \hat\delta_w)=- X^*\Cov(\hat\theta_w, \hat\delta_w)$. Hence, 
\begin{align}\label{eq:adjust}
    \rho^2(\delta_w^*, \hat\delta_w)=\frac{(X^*)^2\Cov^2(\hat\theta_w, \hat\delta_w)}{\left[\Var(Y^*)+ (X^*)^2\Var(\hat\theta_w)\right]\Var(\hat\delta_w)}=\gamma \rho^2(\hat\theta_w, \hat\delta_w).
\end{align}
Consequently, the identity \eqref{eq:corr4} holds for both estimation and prediction, implying the same trade-off between optimal learning and relevant error assessment. 

Section~\ref{sec:mainresult} below will provide a general inequality that captures this trade-off under squared loss, for which identity \eqref{eq:corr4} is a special case. But before presenting that result, we must ask: if we cannot relevantly assess the actual error $\delta$, then what kind of errors have we been assessing? And that is exactly one of the two questions raised in the title of \cite{bates2024cross}: Cross-validation: what does it estimate and how well does it do it? The following section supplements \cite{bates2024cross} to answer this question more broadly and more pedagogically. 

\section{Jay Leno's irony and a statistical magic}\label{sec:Jay}
During one of the years the United States census took place (likely 2000-2001), comedian Jay Leno brought up the issue of under-counting on his \textit{Tonight Show}. He began by informing the audience that the U.S. Census Bureau had just reported that approximately $p$ percentage of the population had not been counted. With an arch smile, he then quipped,  ``But I don't understand---if they knew they missed $p$ percentage of people, why didn't they just add it back?" (The actual value $p$ he used now lies deep in my memory.)

The audience was amused, as was I, though perhaps for different reasons—what amused me was the very appearance of such a nerdy joke on a mainstream comedy show. Humor is often rooted in life's ironies, and whoever crafted this joke clearly understood the irony in announcing both an estimate and its error. In the case of the U.S. Census, the irony—or more accurately, the \textit{magic}—is not as profound as it may seem. The estimation of undercount relies on external data, such as demographic analysis, post-enumeration surveys, administrative records, and other sources. The term \textit{magic} is used here because statistical inference can appear magical to uninitiated yet inquisitive minds. How can one estimate an unknown quantity, and then estimate the error of that estimation, without any external knowledge of the true value?

The magic begins with a sleight of hand—in this case, the word \textit{error} does not  refer to the actual error, as a layperson might assume. Instead, we aim to understand the statistical properties of the actual error by imagining its variations across hypothetical replications. The construction of these replications depends on the philosophical framework one subscribes to, with the two main schools being frequentist and Bayesian (but see \cite{Lin2024To} for a spectrum between them).  Perhaps surprisingly, the key to resolving the apparent paradox in Section~\ref{sec:paradoz} lies in adopting insights from both perspectives.

To see this, consider again the normal example where the true error is \(\delta = \bar{X}_n - \mu\). In the frequentist framework, the hypothetical replications consist of all possible copies of $D={\mathbf X}=\{X_1, \ldots, X_n\}$  generated from \(N(\mu, \sigma^2)\) with the \textit{same} but unknown parameter values \(\theta = \{\mu, \sigma^2\}\). In this replication setting, the expected value of \(\delta^2\), which is the sampling variance of \(\bar{X}_n\), equals \(\sigma^2/n\). It is well-known that under the same replication framework, the expectation of \(\hat{\delta}^2 = S^2_n/n\) is also \(\sigma^2/n\). 

Thus, while \(\delta^2\) and \(\hat{\delta}^2\) are independent of each other for any given \(\theta = \{\mu, \sigma^2\}\), they share the same expectation within the frequentist framework. By invoking the same leap of faith that underpins the frequentist approach—trusting and transferring average behaviors to assess individual cases—we justify \(\hat{\delta}^2\) as an estimate of \(\delta^2\). Such a leap of faith exists regardless of the goal of our data exercise, be it prediction, estimation, or attribution (significance testing), albeit with increased levels of intolerance to the inaccuracy in error assessing, as revealed by the insightful article of \cite{efron2020prediction}.    

For Bayesians, such a leap of faith is unconvincing or even ``irrelevant" in the sense of \cite{Dempster1963c}, as the actual error can differ significantly from its expectation. The independence between \(\hat{\delta}^2\) and \(\delta^2\) suggests that accepting this leap would require a religious level of faith.  In the Bayesian framework, the relevant hypothetical replications include all possible values of \(\theta = \{\mu, \sigma^2\}\) (and their associated probabilities) that could have generated the same data set \(D\), and therefore the same \(\{\bar{X}_n, S^2_n\}\). 

However, for such a replication setting to be realized—for instance, via a simulation—a prior distribution for \(\theta = \{\mu, \sigma^2\}\) must be assumed. This postulation represents the Bayesian leap of faith in actual implementations, since it is virtually certain that a part of the assumption is faith-based  instead of knowledge-driven; for a broader discussion on the necessity of such leaps across all major schools of statistical inference—Bayesian, Fiducial, and Frequentist (BFF)—see \cite{craiu2023six} and more comprehensively the \textit{Handbook on BFF Inference} edited by \cite{berger2024handbook}.

Although we shall not take a Bayesian excursion here, we can borrow the Bayesian concept of allowing \(\theta = \{\mu, \sigma^2\}\) to have a distribution in order to establish a joint replication setting, where both \(D\) and \(\theta = \{\mu, \sigma^2\}\) vary. This framework is relevant (for frequentists) when recommending the same statistical procedure across multiple studies with normal data, where both \(D\) and \(\theta = \{\mu, \sigma^2\}\) may differ from study to study. In the machine learning world—or any domain reliant on training data—such a joint replication setting can be visualized as potential training datasets drawn from related populations, which makes transfer learning a meaningful endeavor \citep[\textit{e.g.},][]{abba2024bayesian}.

For our normal example, given any proper prior on \(\theta\), it can be shown (see Appendix B) that over any proper joint replication of $\{D, \theta\}$, 
\begin{equation}\label{eq:joincorr}
    \rho(\hat\delta^2, \delta^2)=\frac{\gamma^2_{\sigma^2}}{\sqrt{\frac{n+1}{n-1}\gamma^2_{\sigma^2}+\frac{2}{n-1}}\sqrt{3\gamma^2_{\sigma^2}+2}},
\end{equation}
where \(\gamma_{\sigma^2}\) is the coefficient of variation of \(\sigma^2\) with respect to the (proper) prior distribution of \(\sigma^2\).  This correlation is non-negative, providing a plausible measure of how relevant \(\hat{\delta}^2\) is for assessing \(\delta^2\). It is zero if and only if \(\Var(\sigma^2) = 0\), meaning that we revert to the situation of conditioning on a fixed \(\sigma^2\): since \(S_n^2\) is invariant to \(\mu\), \(\hat{\delta}^2\) and \(\delta^2\) remain independent when conditioned on \(\sigma^2\) alone. The fact that \eqref{eq:joincorr} is a monotonic increasing function of \(\gamma_{\sigma^2}\) implies that the relevance of \(\hat{\delta}^2\) for assessing \(\delta^2\) increases as the heterogeneity among the studies—in terms of the within-study variation indexed by \(\sigma^2\)—grows. This monotonicity is intuitive, given that \(S_n^2\) is an unbiased and asymptotically efficient estimator of \(\sigma^2\), and \(\hat{\delta}^2\) is useful for comparing the magnitudes of \(\delta^2\) across studies with different \(\sigma^2\) values. However, the fact that this correlation can never exceed \(1/\sqrt{3} \approx 0.577\) is unexpected. For those of us who believe that mathematical results are never coincidental, contemplating the intricacies of this bound might induce insomnia (while serving as a cure for many others).

This joint replication framework clarifies the role of \(\hat{\delta}^2\) as an \textit{adaptive benchmark} for assessing the statistical properties of \(\delta^2\) over the hypothetical replications. \textit{That} is statistical magic—the ability to establish cross-study comparisons based on a single study. More broadly, the magic lies in creating hypothetical ``control'' replications \(\{\tilde{D}, \tilde{\theta}\}\) from the actual ``treatment'' \(\{D, \theta\}\) at hand, as elaborated in \cite{liu2016there}, borrowing the metaphor of individualized treatment.

Generally speaking, the magic relies on two tricks: (I) creating replications within \(D\), and (II) linking those replications to the imagined variations of \(D\) through the within-\(D\) replications from (I). The first trick is applicable when the mechanism generating the data \(D\) inherently includes (higher resolution) replications, either by design (\textit{e.g.}, simple random sampling) or by declaration (\textit{e.g.}, imposing an i.i.d. structure as a working assumption). The second trick is enabled by theoretical understanding (\textit{e.g.}, the relationship between the distribution of the sample mean and the distribution of the individual samples) or by simulations and approximations that are enabled by (I), such as the Bootstrap \citep[see][for a discussion]{craiu2023six}.

The magic metaphor also serves as a reminder that magic relies on illusions, and interpreting average errors as actual ones is such an illusion. With that understanding, we might wonder if it’s possible to assess the actual error with greater relevance. For example, in the normal case, one might ask whether a different error estimate \(\check{\delta}\) could be more relevant for \(\delta = \bar{X}_n - \mu\), in the sense that \(\rho(\check{\delta}, \delta) > 0\) given any value of $\theta=\{\mu, \sigma^2\}$. The classical statistical literature offers a fairly clear answer to this question, as discussed below.

\section{From UMVUE to an uncertainty principle for unbiased learning}\label{sec:mainresult}

The celebrated Cramér–Rao bound, more broadly known as the information inequality \citep[see][Ch. 2]{lehmann2006theory}, tells us that if \(\hat{\theta}\) is an unbiased estimator for \(\theta\) under a parametric model \(f(D|\theta)\), then under mild conditions, \(\Var(\hat{\theta}) \geq I^{-1}(\theta)\), where \(I(\theta)\) is the expected Fisher information. For the normal example, when we take \(\theta = \mu\) (temporarily assuming \(\sigma^2\) is known), we have \(\Var(\bar{X}_n) = \sigma^2/n = I^{-1}(\mu)\), where \(I(\mu)\) is the expected Fisher information from \(f(X_1, \ldots, X_n | \mu)\). Thus, we know \(\bar{X}_n\) is UMVUE for \(\mu\).

It is well-known that an estimator \(\hat{\theta}\) is UMVUE if and only if it is uncorrelated with any unbiased estimator \(U\) for zero for any \(\theta\) \citep[see][Ch. 2]{lehmann2006theory}, that is, \(\E_\theta[(\hat{\theta} - \theta) U] = 0\), whenever $\E_\theta(U)=0$. Since \(\hat{\theta} - \theta\) is simply the actual error \(\delta\), this result implies that conditioning on \(\theta\), it is impossible to have an error assessment \(\hat{\delta}\) for \(\delta\) that is both unbiased and relevant at the same time, \textit{i.e.}, \(\E_\theta(\hat{\delta}) = 0\) and \(\rho_\theta^2(\hat{\delta}, \delta) > 0\) cannot hold simultaneously for any \(\theta\), where we inject the subscript $\theta$ in $\rho_\theta$ to explicate that the correlation is with respect to $f(D|\theta)$ for fixed $\theta$.

Intuitively, if any unbiased error assessment \(\hat{\delta}\) is correlated with \(\delta\), then some part of the actual error \(\delta\) is predictable by \(\hat{\delta}\). This means that we could improve \(\hat{\theta}\) without losing its unbiasedness, which contradicts the fact that \(\hat{\theta}\) is already an UMVUE. An astute reader may quickly recognize that this insight has much broader implications than merely for UMVUEs. The following result is a proof of this realization,  using the same  proof strategy as for UMVUE, but establishes a broader quantitative result than the aforementioned qualitative ``if and only if" result for UMVUE. The result is presented in the scalar case for simplicity, but its multivariate counterpart can be derived easily using corresponding matrix notation.

Specifically, let \(Q \in \mathbb{R}\) be our target of learning, which could represent a future outcome, a model parameter, a latent trait, etc. Suppose the state space of our data \(D\) is \(\Omega\) and \(\hat{Q}: \Omega \rightarrow \mathbb{R}\) is our learning algorithm, or a \textit{learner}  for \(Q\). For any learner \(\hat{Q}\), let \(\hat{\delta}_{\hat{Q}}: \Omega \rightarrow \mathbb{R}\) be an assessment (\textit{e.g.}, an estimator) of the exact (additive) error of \(\hat{Q}\), namely, \(\delta_{\hat{Q}} = \hat{Q} - Q\).  Let \(L(\hat{Q}, Q)\) be the loss function, and \({\cal P} = \{P_s(D; Q), s \in S\}\) be  the family of distributions under which we calculate the learning risk: \(R_s(\hat{Q}) = \E_s[L(\hat{Q}, Q)]\). Note that \(Q\) may be a function of \(s\) (\textit{e.g.}, when estimating the model parameter \(s\)) or it may be a random variable itself (\textit{e.g.}, a future realization), in which case the notation \(P_s(D; Q)\) represents the joint distribution over \(D\) and \(Q\).

\begin{theorem}\label{th:main}
Let \(L(\hat{Q}, Q) = (\hat{Q} - Q)^2\) be the squared loss, and let \(L^2_{\cal P}\) denote the collection of all square-integrable functions with respect to \({\cal P}\). Define
\begin{equation}\label{eq:defineQ}
{\cal Q} = \{\hat{Q}(D) \in L^2_{\cal P} :  \E_s(\hat{Q} - Q) = 0, \forall s \in S\}
\end{equation}
as the collection of unbiased learners of \(Q\) with respect to \({\cal P}\). For any \(\hat{Q} \in {\cal Q}\), define
\begin{equation}\label{eq:definedelta}
{\cal E}(\hat{Q}) = \{\hat{\delta}_{\hat{Q}}(D) \in L^2_{\cal P} :  \E_s(\hat{\delta}_{\hat{Q}}) = 0, \forall s \in S\}
\end{equation}
as the collection of corresponding unbiased error assessors for \(\delta_{\hat{Q}}\). Suppose there exists an optimal learner \(\hat{Q}^{\rm opt} \in {\cal Q}\), with risk \(R_s^{\rm opt} < \infty\) under \(f_s, s \in S\). Then:
\begin{description}
    \item[(I)] For any \(\hat{Q} \in {\cal Q}\) and any corresponding \(\hat{\delta}_{\hat{Q}} \in {\cal E}(\hat{Q})\), we have
    \begin{equation}\label{eq:corrinq}
    \rho_s^2(\delta_{\hat{Q}}, \hat{\delta}_{\hat{Q}}) \leq \frac{R_s(\hat{Q}) - R_s^{\rm opt}}{R_s(\hat{Q})} \equiv RR_s(\hat{Q}), \quad \forall s \in S,
\end{equation}
    where \(RR_s(\hat{Q})\) is the relative regret of \(\hat{Q}\) under distribution \(P_s\), and it is set to zero if \(R_s(\hat{Q}) = 0\).
    \item[(II)] Equality \(\rho_s^2(\delta_{\hat{Q}}, \hat{\delta}_{\hat{Q}}) = RR_s(\hat{Q})\) holds for any particular \(s\) if and only if \(R_s^{\rm opt}\) is attainable in the sub-class \({\cal Q}(\hat{Q}, \hat{\delta}_{\hat{Q}}) = \{\hat{Q} - \lambda \hat{\delta}_{\hat{Q}} :  \forall \lambda \in \mathbb{R}\} \subset {\cal Q}\).
\end{description}
\end{theorem}

\begin{proof}
For any given \(\hat{Q} \in {\cal Q}\) (which is non-empty since \(\hat Q^{\rm opt} \in {\cal Q}\)) and any \(\hat{\delta}_{\hat{Q}} \in {\cal E}(\hat{Q})\) (which is non-empty since \(\hat{\delta}_{\hat{Q}} \equiv 0\) is always included), we define \(\hat{Q}_\lambda = \hat{Q} - \lambda \hat{\delta}_{\hat{Q}}\) for any constant \(\lambda \in \mathbb{R}\). Under our assumptions, \(\E_s(\hat{Q}_\lambda - Q) = 0\), and \(\hat{Q}_\lambda \in L^2_{\cal P}\),  implying  \(\hat{Q}_\lambda \in {\cal Q}\). Since \(\hat{Q}_\lambda - Q = \delta_{\hat{Q}} - \lambda \hat{\delta}_{\hat{Q}}\) and it has mean zero under \(f_s(D; Q)\), we have
\begin{equation}\label{eq:proofine1}
R_s^{\rm opt} \leq R_s(\hat{Q}_\lambda) = \Var_s(\delta_{\hat{Q}} - \lambda \hat{\delta}_{\hat{Q}}) = \Var_s(\delta_{\hat{Q}}) + \lambda^2 \Var_s(\hat{\delta}_{\hat{Q}}) - 2 \lambda \Cov_s(\delta_{\hat{Q}}, \hat{\delta}_{\hat{Q}}), \quad \forall s \in S.
\end{equation}
Since the left-hand side of this inequality is free of \(\lambda\), the inequality holds when we minimize the right-hand side over \(\lambda \in \mathbb{R}\), which is achieved at \(\lambda = \lambda^* = \Cov_s(\delta_{\hat{Q}}, \hat{\delta}_{\hat{Q}})/\Var_s(\hat{\delta}_{\hat{Q}})\), assuming \(\Var_s(\hat{\delta}_{\hat{Q}}) > 0\). (When \(\Var_s(\hat{\delta}_{\hat{Q}}) = 0\),  \(\rho_s(\delta_{\hat{Q}}, \hat{\delta}_{\hat{Q}}) = 0\); hence \eqref{eq:corrinq} holds trivially, and we can set \(\lambda^* = 0\).) Thus, we obtain
\[
R_s^{\rm opt} \leq \Var_s(\delta_{\hat{Q}}) \left[1 - \rho^2_s(\delta_{\hat{Q}}, \hat{\delta}_{\hat{Q}})\right], \quad \forall s \in S,
\]
which yields  \eqref{eq:corrinq} since \(R_s(\hat{Q}) = \Var_s(\delta_{\hat{Q}})\) when \(\E_s(\delta_{\hat{Q}}) = 0\). This proves part (I).

Part (II) follows from \eqref{eq:proofine1} as well, because the equality holds  there if and only if \(R_s^{\rm opt}\) is attainable by \(\hat{Q}_{\lambda^*} \in {\cal Q}(\hat{Q}, \hat{\delta}_{\hat{Q}})\). This includes the case with $\Var_s(\hat{\delta}_{\hat{Q}}) = 0$, where the result holds trivially, because then  \(\rho_s(\delta_{\hat{Q}}, \hat{\delta}_{\hat{Q}}) = 0\) and \(R_s(\hat{Q}) = R_s^{\rm opt}\), \textit{i.e.}, \(\hat{Q}\) itself is optimal.
\end{proof}

The immediate implication of inequality \eqref{eq:corrinq} is that there is no free lunch. If we want to increase the relevance of our assessment \(\hat{\delta}_{\hat{Q}}\) for the actual error \(\delta_{\hat{Q}}\) by increasing their correlation, we must also increase the relative regret for \(\hat{Q}\), effectively sacrificing degrees of freedom of learning for the error assessment. Conversely, the less regret in \(\hat{Q}\), the less relevant its error assessment will be to the actual error. In the extreme case, when \(\hat{Q} = \hat{Q}^{\rm opt}\), we arrive at the following result, where by a \textit{relevant error assessor} we mean it is linearly correlated with the actual error of the learner. 

\begin{corollary}
Under the same setup as in Theorem~\ref{th:main}, the following two assertions cannot hold simultaneously:
\begin{itemize}
    \item[(A)] \(\hat{Q} \in {\cal Q}\) is an \textit{optimal and unbiased learner} for \(Q\) under \(P_s\); and
    \item[(B)] \(\hat{Q}\) has an \textit{unbiased and relevant error assessor }\(\hat{\delta}_{\hat{Q}} \in {\cal E}(\hat{Q})\).
\end{itemize}
\end{corollary}

\section{Beyond unbiased learning and error assessing}\label{sec:asymresults}
A key limitation of Theorem~\ref{th:main} is the requirement that both the learner and error assessor must be unbiased. An immediate generalization is to consider cases where both are asymptotically unbiased, under an asymptotic regime with respect to some information index $\iota$, such as the size of data. Mathematically, given a sequence of error order $e_\iota$ such that $\limsup_{\iota\rightarrow\infty} |e_\iota|=0$, we can modify the classes of the learners and error assessors in \eqref{eq:defineQ} and \eqref{eq:definedelta} respectively by 
\begin{align}\label{eq:defineQe}
{\cal Q}_\iota=&\{\hat Q(D) \in L^2_{\cal P} :  \E_s[\hat Q(D)-Q]=O(e_\iota), \forall s\in S\},\\
\label{eq:definedeltae}
{\cal E}_{\iota}(\hat Q)=&\{\hat \delta_{\hat Q}(D) \in L^2_{\cal P}:  \E_s(\hat \delta_{\hat Q})=O(e_\iota), \forall s\in S\},
\end{align}
where $O(e_\iota)$ is the standard notation of being on the same order as $e_\iota$. That the error assessor $\hat\delta_{\hat Q}$ must share the same order of expectation as the actual error $\delta_{\hat Q}$ is a necessary requirement to render the term `error assessor' meaningful, as otherwise anything could be regarded as $\hat\delta_{\hat Q}$.  With these modifications, we have the following asymptotic counterpart of Theorem~\ref{th:main}. 
\begin{theorem}\label{th:maina} Assume the same setup as Theorem~\ref{th:main}, but with 
${\cal Q}$ and ${\cal E}(\hat Q)$ extended respectively to ${\cal Q}_\iota$ and ${\cal E}_{\iota}(\hat Q)$.  We then have  
\begin{equation}\label{eq:corrinqa}
\rho^2_s(\delta_{\hat Q}, \hat \delta_{\hat Q}) \le RR_s(\hat Q) + O(e^2_\iota), \quad \forall s\in S,
\end{equation}
where $e_\iota$ is a sequence of vanishing error rates that determines the asymptotic regime. 
\end{theorem}
\begin{proof} For $\hat Q \in {\cal Q}_\iota$ and $\hat \delta_{\hat Q} \in {\cal E}_\iota(\hat Q_\iota)$, we can write $\E_s(\delta_{\hat Q})=a_\iota$ 
and $\E_s(\hat\delta_{\hat Q})=b_\iota$ where $a_\iota=O(e_\iota)$ and $b_\iota=O(e_\iota)$ by our assumption. Hence for $\hat Q_\lambda=\hat Q - \lambda\hat\delta_{\hat Q}$, $\E_s(\hat Q_\lambda-Q)=a_\iota - \lambda b_\iota=Q(e_\iota)$ for any $\lambda$, implying that $\hat Q_\lambda\in {\cal Q}_\iota. $  Let $\lambda^*$ be the minimizer of $\Var_s\left[\delta_{\hat Q} - \lambda\hat\delta_{\hat Q}\right]$, as defined in the proof of Theorem~\ref{th:main}. The optimality of $R^{\rm opt}_s$ then implies that
\begin{align*}
R^{\rm opt}_s\le R_s(\hat Q_{\lambda^*})&= \Var_s\left[\delta_{\hat Q} - \lambda^*\hat\delta_{\hat Q}\right]
+ \left[\E_s(\delta_{\hat Q} - \lambda^*\hat\delta_{\hat Q})\right]^2 \\
&=\Var_s(\delta_{\hat Q})\left[1 - \rho^2_s(\delta_{\hat Q}, \hat \delta_{\hat Q}) \right]
+ (a_\iota-\lambda^* b_\iota)^2.\\
&\le R_s(\delta_{\hat Q})\left[1 - \rho^2_s(\delta_{\hat Q}, \hat \delta_{\hat Q}) \right]
+  (a_\iota-\lambda^* b_\iota)^2.
\end{align*}
But this proves the inequality \eqref{eq:corrinqa}  because $(a_\iota-\lambda^* b_\iota)^2=O^2(e_\iota)= O(e_\iota^2)$. 
\end{proof}

A major application of Theorem~\ref{th:maina} is for the maximum likelihood estimator $\hat Q_{\rm MLE}$, which under regularity conditions is efficient and asymptotically normal \citep[\textit{e.g.},][]{lehmann2006theory} and hence it is asymptotically optimal under the squared loss. Theorem~\ref{th:maina} says that asymptotically, there cannot be any relevant error assessor $\hat\delta_{\rm MLE} \in {\cal E}_\iota(\hat Q_{\rm MLE})$  that is asymptotically correlated with the actual error $\delta_{\rm MLE}=\hat Q_{\rm MLE} - Q$.  When $\{\hat\delta_{\rm MLE}, \delta_{\rm MLE}\}$ are jointly asymptotically normal, then Theorem~\ref{th:maina} would imply that  any such $\hat\delta_{\rm MLE}$  will be asymptotically independent of $\delta_{\rm MLE}$.  It is worthy noting that the same would hold for any estimator that is asymptotically normal and optimal (under quadratic loss), such as those studied in the classic work by \cite{wald1943tests} and \cite{le1956asymptotic}.

Because the asymptotic variance of the MLE can be well approximated by the inverse of Fisher information, especially the observed Fisher information \citep{efron1978assessing},  the preceding result might lead some readers to wonder if the MLE and the observed Fisher information are asymptotically independent, or at least the MLE and the inverse of the observed Fisher information $I_{\text obs}^{-1}(\hat Q)$  are asymptotically uncorrelated. The normal example given in Section~\ref{sec:paradoz} may be especially suggestive,  since the MLE for $\mu$, $\bar X_n$,  is independent of   $I_{\text obs}^{-1}(\hat\mu)=n/\hat \sigma^2_{\rm MLE}=n^2/[(n-1)S^2_n]$.   However it will be a mistake to generalize from this example.  

Consider the same normal model $N(\mu, \sigma^2)$,  but our goal now is  to estimate  the variance $\sigma^2$.  The MLE for $\sigma^2$ is $\hat\sigma^2_{\rm MLE}=(n-1)S^2_n/n$, and the  corresponding observed Fisher information (pretending $\mu$ is known) is $I^{-1}(\hat\sigma^2_{\rm MLE})=2\hat\sigma^4_{\rm MLE}/n$; hence they have a deterministic relationship.  However, this is not a contradiction to Theorem~\ref{th:maina} because $I^{-1}(\hat\sigma^2_{\rm MLE})$ is not an unbiased assessment of the actual error, but rather its variance.  Since the variance is effectively an index of the \textit{problem difficulty} for estimation \citep[as termed in][]{meng2018statistical}, it is entirely natural to expect that the variance can vary closely with the value of the estimand.  The normal mean problem is a special case because it is a location family, for which shifting the mean only changes the value of  the estimand, but does not alter the difficulty of its estimation.  This point is reinforced  if we reparameterize $\sigma^2$ via  $\eta=\log\sigma^2$, which yields $\hat\eta_{\rm MLE}=\log \hat\sigma^2_{\rm MLE}$ and   $I^{-1}(\hat\eta^2_{\rm MLE})=2/n$, and they are now trivially independent of each other,  because $\hat\eta_{\rm MLE} - \eta \sim \log \chi^2_{n-1}-\log{(n-1)} $ is a location family.  

The consideration of the relationship between the MLE and the Fisher information provides a natural segue to the following discussion involving the relationship between inequity \eqref{eq:corrinq} and the Cram\'er-Rao low bound. As is well documented\footnote{See the video on C.R. Rao: A Life in Statistics II at \text{\tt https://www.youtube.com/watch?v=eaxjUxoCx5w\&t=324s}}, the seminal work by \cite{rao1945information} was prompted by a  question raised during a lecture Rao gave in 1943 on whether there could be a small-sample counterpart of the asymptotic efficiency for MLE as captured by the Fisher information.  However, the significance of this work goes beyond accenting the role of Fisher information, because the Cram\'er-Rao inequality can be viewed as a statistical counterpart of the fundamental Heisenberg Uncertainty Principle (HUP, \cite{griffiths2018introduction}) via the notion of \textit{co-variation}, as explored in the next three sections. 

\section{Measuring co-variation without probabilistic joint-state specifications}\label{sec:definingcov}

In statistical and (ordinary) probabilistic literature, the most commonly adopted measure of the co-variation of two real-valued random variables $G$ and $H$ is their covariance $\Cov(G, H)$ (which includes correlation once $G$ and $H$ are standardized) defined via their joint probabilistic distribution $F_{G, H}(g, h)$:
\begin{equation}\label{eq:probcov}
\Cov(G, H)=\int\int (g-\mu_G)(h- \mu_h)F_{G, H}(dg, dh)
=\langle (g-\mu_G), (h- \mu_h) \rangle_F,
\end{equation}
where $\mu_G$ and $\mu_H$ are respectively the means of $G$ and $H$,
which, without loss of generality, we will assume to be zero for the subsequent discussions for notational simplicity. The subscript $F$ in the inner product notation highlights the critical dependence of $\Cov(G, H)$ on their joint distribution $F(g, h)$.  The elegant Hoeffding identity \citep{Hoeffding1940} 
\begin{equation}\label{eq:hoeffding}
\Cov(G, H)=\int\int \left[F_{G. H}(g, h) - F_{G}(g)F_H(h)\right] dg dh, 
\end{equation}
where $F_G$ and $F_H$ are the marginal (cumulative) distributions, further highlights how the covariance measures the co-variation in
$G$ and $H$ as captured by their joint distribution, with respect to their benchmark distribution under the assumption of independence. 

For HUP, it seems natural to take $G=x$, the position of a particle,  and $H=p$, its momentum, to follow the standard notation in quantum mechanics. It is textbook knowledge \citep[\textit{e.g.}][]{landau2013quantum,griffiths2018introduction} that densities of the position $x$ and momentum $p$ are given by $|\psi(x)|^2$ and $|\varphi(p)|^2$ respectively,
where $\psi(x)$ is a complex valued position wave function, and the momentum wave function $\varphi(p)$ is a scaled Fourier transform of $\psi(x)$ in the form of 
\begin{equation}\label{eq:Fourier}
    \varphi(p)= \frac{1}{\sqrt{2\pi \hbar}} \int_{-\infty}^{\infty} \psi(x) \, e^{-ipx / \hbar}dx,
\end{equation}
where the scale factor $\hbar=h/(2\pi)$, with $h=6.6260701\times 10^{-34}$, the Planck's constant. Clearly, $\psi(x)$ is the inverse Fourier transform of $\varphi(p)$, and together $x$ and $p$ form a pair of the so-called conjugate variables \citep{stam1959some}. 

As a statistician, once I understood how the marginal distributions for $x$ and $p$ were constructed, I naturally asked for their joint distribution.  This is where things become intriguing or puzzling to those of us who are trained to model non-deterministic relationships via probability, because (quantum) physicists' answer would be that there is no joint probability distribution for $x$ and $p$---not that they are unknown, but that there cannot be one. Unlike the mystery of deep learning to statisticians---and its winning of the Nobel prize in physics only makes it more intriguing or puzzling---I found good clues to the inadequacy of ordinary probability for dealing the quantum world by the very fact that its mathematical modeling involves non-commutative relationships, such as between operators or matrices.  

Perhaps the easiest way to see potential complications with non-commutative relationship is to consider the problem of generalizing the notion of variance to co-variance with complex-valued variables. With real-valued random variables $G$ and $H$ having a joint distribution $F$, we know variance is the co-variance of a variable with itself, that is,  $\Var(G)=\Cov(G, G)$. In other words, when we link variance with an inner product, \textit{i.e.}, $\Var(G)=\langle G , G \rangle_F$, there is a natural extension for covariance by defining $\Cov(G, H)=\langle G , H\rangle_F$. However, with the ordinary definition of the co-variance, this extension works only if the inner product is symmetric, that is, $\langle G , H\rangle_F=\langle H, G\rangle_F$, since $\Cov(G, H)=\Cov(H, G)$ in the \textit{real} world.

This is where the complex world is, literally, more complex than the real world. For two complex-valued $L^2$ functions $u(y)$ and $v(y)$ on $y\in \Omega$, the inner product is not symmetric, because it is defined by
\begin{equation}\label{eq:covariation}
\langle u | v \rangle_\mu \equiv  \int_\Omega \bar u(y)v(y) \mu(dy) \not = \langle v | u\rangle_\mu \equiv  \int_\Omega \bar v(y)u(y) \mu(dy),
\end{equation}
where $\bar u$ is the complex conjugate of $u$, and $\mu$ is a baseline measure,  which does not need to be a probabilistic measure.  This non-commutative property is at the heart of quantum mechanics, as reviewed in the next Section. It can also been seen with matrix mechanics,  since for any two matrices $A$ and $B$ or more broadly operators, in general $AB\not=BA$.  The very fact that a regular joint probability specificity must render $\Cov(u, v)=\Cov(v, u)$ should remind us that whatever `joint specification' of $u$ and $v$ we come up with, it will be more nuanced than a direct probabilistic distribution for $\{u, v\}$ whenever \eqref{eq:covariation} rears its head. This phenomena is not unique to the quantum world, since a similar situation happens with the notion of quasi-score functions, which can violate a symmetry requirement for genuine score functions, as reviewed in Appendix C. 

However, this complication does not imply that probabilistic thinking is out the window. Because $\overline{\langle v | u \rangle_\mu} =\langle u| v \rangle_\mu $, we see that if we define
${\Cov}(u, v)=\langle u| v \rangle_\mu$, then its magnitude,
$|\Cov(u, v)|=|\Cov(v, u)|$  
is symmetric. Therefore, as long as $|\Cov(u, v)|$ is used as a measure of the magnitude of the co-variation between $u$ and $v$, we can treat it as if it were the magnitude of a standard  probabilistic co-variance.  In other words, the concept or at least the essence of \textit{co-variance} can be extended to non-probabilistic settings, and this extension perhaps can help our appreciation of HUP from a statistical perspective, as detailed in the next Section.

\section{A lower resolution co-variation: co-variance of generating mechanisms}\label{sec:covoperator}

In the quantum world, we have seen that a particle's position and momentum have their respectively well-defined probability distribution,  and we can express $\Var(x)={\langle f| f \rangle_\mu}$ and $\Var(p)={\langle g|g\rangle_\mu}$, where $f(x)=x\psi(x)$ and $g(p)=p\phi(p)$. It is then mathematically tempting to define $\Cov(x, p)={\langle f| g \rangle_\mu}$ and $\Cov(p, x)={\langle g| f\rangle_\mu}$, using the notation of the pervious section. This construction is problematic starting from the very notation $\Cov(x, p)$, since it may suggest that we are measuring the co-variance between the position and momentum as \textit{states}, which creates an epistemic disconnect with the understanding that \textit{a joint statehood} of $x$ and $p$ does not exist or cannot be constructed in the quantum world.

However, $x$ and $p$ clearly have physical relationships. Indeed the so-called Stam's uncertainty principle \citep{stam1959some} establishes that 
\begin{equation}\label{eq:stam}
C^2\Var(x) - J(p) \ge 0 \quad {\rm and} \quad C^2\Var(p) - J(x) \ge 0,
\end{equation}
where $C=4\pi$ for standard Fourier transform, and $C=2/\hbar$ when we use the $\hbar$-scaled Fourier transform \eqref{eq:Fourier}. Here $J(p)$ is the Fisher information for the density of $p$, $f(p)$, that is, 
\begin{equation}\label{eq:Fisher}
    J(p) = \int_{-\infty}^\infty \left[\frac{d\log f(p)}{dp}\right]^2 f(p) dp,
\end{equation}
and similarly for $J(x)$. For readers who are unfamiliar with defining Fisher information for a density itself instead of its parameter, 
$J(p)$ is the same as the Fisher information for the location family $f(p-\theta)$, where $\theta$ shares the same state space as $p$ (in the current case, the real line). In the same vein, the Cram\'er-Rao inequality can be applied to the density itself, which leads to $\Var(x)\ge J^{-1}(x)$ and $\Var(p)\ge J^{-1}(p)$. Consequently,
as shown in \cite{dembo1990information} and \cite{dembo1991information}, 
\begin{equation}\label{eq:HUP}
\Var(x)\Var(p)\ge C^{-2}=\frac{\hbar^2}{4},
\end{equation}
which is the same as the usual expression of HUP proved in \cite{kennard1927}:  
\begin{equation}\label{eq:keen}
    \Delta x \Delta p \ge \frac{\hbar}{2},
\end{equation}
where $\Delta x$ and $\Delta p$ denote respectively the standard deviation of $x$ and $p$. \cite{dembo1990information} and \cite{dembo1991information} also used \eqref{eq:stam} to prove that HUP implies the Cram\'er-Rao inequality.

The Stam's uncertainty principle is elegant, and it reveals a kind of relationship between two marginal distributions that is not commonly studied in statistical literature, because it bypasses the specification of a joint distribution between $x$ and $p$. However, this does not rule out---and indeed it suggests---that we can consider quantifying the relationships between the \textit{mechanisms} that generate $x$ and $p$. A mechanism can generate a single state, many states, or no states at all---which is equivalent to presenting itself as a whole---at any given circumstance, such a temporal instance. Hence quantifying relationships among mechanisms is a broader construct than that for the states they generate. 

For statistical readers, a reasonable analogy is to think about the notion of \textit{likelihood}. When we employ a likelihood, we can consider a single likelihood value (\textit{e.g.}, at the MLE), several likelihood values (\textit{e.g.}, likelihood ratio tests), or not any particular value but the likelihood function as a whole (\textit{e.g.}, for Bayesian inference). By considering co-variations at the (resolution) level of mechanisms instead of states, we may find it less foreign to contemplate indeterminacy of relationship, such as between two sets---including empty ones---of the states generated by related mechanisms. 

Of course, one may wonder if any relationship between two mechanisms itself can be indeterminable. The logical answer is yes, but fortunately for quantum mechanics we do not need go that far. As any useful quantum mechanics textbook \citep{landau2013quantum,griffiths2018introduction} teaches us, the position mechanism and momentum mechanism can be  represented mathematically via the so-called position operator $\hat x$ and momentum operator 
$\hat p$, to follow the notation in quantum mechanics, and they are tethered together when being applied to the same wave function $\psi(x)$ (in the position space\footnote{One can define the operators equivalently in the conjugate momentum space via $\hat p\circ \varphi(p) = p\varphi(p)$ and $\hat x \circ \varphi(p) = i \hbar \varphi^\prime(p)$, where the momentum wave function $\varphi(p)$ is the Fourier transform of $\psi(x)$ \citep{griffiths2018introduction}.}), that is 
\begin{equation}\label{eq:operator}
    \hat x\circ \psi(x) = x \psi(x), \quad {\rm and} \quad  \hat p\circ \psi(x) = - i \hbar \psi^\prime(x).
\end{equation}
That is, the position operator acts on $\psi$ by multiplying  $\psi$ with its argument, and the momentum operator acts on $\psi$ by differentiating it, and multiplying it by $-i\hbar$, where $i=\sqrt{-1}$.  

With these representations of the mechanisms, we can measure their co-variations induced by changing the state $x$ in real line (as a univariate case) via the inner products,  with respect to a common measure $\mu$, typically Lebesgue measure. That is, we can define
\begin{align}\label{eq:covcommon}
    \Cov(\hat x, \hat p) &= \langle \hat x\circ \psi | \hat p\circ \psi\rangle_{\mu} = -i\hbar\int_{-\infty}^{\infty} x \bar\psi(x)  \psi^{\prime}(x) \, dx; \\
    \Cov(\hat p, \hat x)&= \overline{\Cov(\hat x, \hat p)}
    =  i\hbar\int_{-\infty}^{\infty} x \psi(x)\bar\psi^{\prime}(x)  \, dx 
    = -i\hbar \left( 1 + \int_{-\infty}^{\infty} x \bar{\psi}(x) \psi^{\prime}(x) \, dx \right). \label{eq:covcommonc}
\end{align}
Here the last equality is obtained by integration by parts and by using the fact that $|\psi(x)|^2$ is a probability density and that $x|\psi(x)|^2$ vanishes at $x=\pm \infty$ (because physicists assume the mean position is finite).  Together, expressions \eqref{eq:covcommon}-\eqref{eq:covcommonc} imply that
\begin{equation}\label{eq:nonc}
    \Cov(\hat x, \hat p) - 
    \Cov(\hat p, \hat x)
    = i\hbar, 
\end{equation}
which is also the consequence of the so-called \textit{canonical commutation relation} \citep{griffiths2018introduction},  
\begin{equation}\label{eq:canonical}
\hat x\circ \hat p - \hat p\circ \hat x =i\hbar,
\end{equation}
which holds because $\hat x\circ (\hat p\circ f(x)) - \hat p\circ (\hat x\circ f(x))=i\hbar f(x)$ for any differentiable function $f$. 

An immediate consequence of \eqref{eq:nonc} is that the magnitude of the covariances between $\hat x$ and $\hat p$ is bounded below regardless of the form of the wave function $\psi(x)$.
This is because for any complex number $z$, $|z|^2\ge |{\rm Im}(z)|^2=|(z - \bar z)/2i|^2$. Hence the identity \eqref{eq:nonc} implies that
\begin{equation}\label{eq:hbound}
|\Cov(\hat x, \hat p)|^2 \ge \left[\frac{ \Cov(\hat x, \hat p) - 
    \Cov(\hat p, \hat x)}{2i}\right]^2=\frac{\hbar^2}{4}.
\end{equation}
As reviewed in the next section, inequality \eqref{eq:hbound} implies HUP in the form of \eqref{eq:keen}, just as Stam's uncertainty principle does.  For that purpose, it is worth pointing out that marginally, 
\begin{align}\label{eq:marginV}
    \Var(\hat x) &= \langle \hat x\circ \psi |  \hat x\circ \psi\rangle_{\mu} = \int_{-\infty}^{\infty} x^2 \bar\psi(x)\psi(x)dx=\int_{-\infty}^{\infty} x^2 |\psi(x)|^2dx; \\
    \Var(\hat p) &= \langle \hat p\circ \psi |  \hat p\circ \psi\rangle_{\mu}
    = \hbar^2\int_{-\infty}^{\infty} \bar\psi^{\prime}(x)\psi^{\prime}(x)  \, dx = \int_{-\infty}^{\infty} p^2|\varphi(p)|^2 dp, \label{eq:marginVp}
\end{align}
where the last equation in \eqref{eq:marginVp} is due to the fact that $\varphi(p)$ is the ($\hbar$-scaled) Fourier transformation of $\psi(x)$, as given in \eqref{eq:Fourier}.
These two equalities tell us that when we consider either the position or the momentum by itself,  its mechanism-level variance, $\Var(\hat x)$ or $\Var(\hat p)$,  and the state-level variance, $\Var( x)$ or $\Var(p)$,  is the same. This renders the unity between the mechanism-level representation (as a distribution or operator) and the state-level representation (as a observable or latent variable), a distinction seldom made conceptually under the ordinary probability framework. However, this distinction can  be crucial once we go outside the regular probability framework, as in the current context of measuring co-variations between the position and momentum.   

\section{Bounding co-variations: A commonality of uncertainty  principles}\label{sec:covariation}

With co-variances constructed broadly, we can study the similarities and differences between inequality \eqref{eq:corrinq} and the Cram\'er-Rao inequality, as well as their intrinsic connections with HUP. 
Specifically, both inequalities are based on bounding joint variations of two random objects, say, $G$ and $H$, by their marginal variations. For \eqref{eq:corrinq},  under the unbiasedness  assumptions and using the notation given in Section~\ref{sec:mainresult}, if we write $G=\delta_{\hat Q}$ and $H=\hat\delta_{\hat Q}$, then inequality \eqref{eq:corrinq} is the consequence of (omitting subscript $s$):
\begin{equation}\label{eq:GH}
   \Cov^2(G, H) \le \Var(H) \left[\Var(G) - R^{\rm opt}\right]. 
\end{equation}
For the Cram\'er-Rao inequality, we can take the same $G=\delta_{\hat Q}=\hat Q - Q$, where $\hat Q$ is an unbiased estimator for $Q$. We then let $H=S(\theta|D)$, the score function from a sampling model of our data $D$, $f(D|\theta)$, with $Q=Q(\theta)$. It is known that the Cram\'er-Rao inequality is the same as \citep[\textit{e.g.},][]{lehmann2006theory}
\begin{equation}\label{eq:Rao}
   [Q^{\prime}(\theta)]^2=\Cov^2(G, H) \le \Var(H) \Var(G),
\end{equation}
where $Q^{\prime}(\theta)$ is the derivative for $Q(\theta)$. (When $Q(\theta)$ is not differentiable, we can apply the bound given by \cite{chapman1951minimum}) in terms of likelihood ratio or elasticity.)

Evidently, inequality \eqref{eq:Rao} is an application of the Cauchy-Schwartz inequality. In contrast, inequality \eqref{eq:GH} delivers a more precise bound because of the subtraction of the term $R^{\rm opt}$. Indeed, inequality \eqref{eq:GH} is often an equality because the condition in (II) of Theorem~\ref{th:main} frequently holds in practice. 
Give the two inequalities share the same type of $G$, the difference must be attributable to something distinctive between the two $H$'s. Whereas both $H$'s have zero expectation, the first $H=\hat\delta_{\hat Q}$ is a \textit{statistic}, required to be a function of data $D$ only. In contrast, the second $H=S(\theta|D)$ is a \textit{random function}, depending on both data $D$ and the unknown $\theta$. Since the actual error $\delta_{\hat Q}=\hat Q - Q(\theta)$ is also a random function, the second $H$ can co-variate with $G$ to a greater extent than the first $H$ can. Consequently, $\Cov^2(G, H)$ can reach a looser upper bound in \eqref{eq:Rao} than in \eqref{eq:GH}.  As an illustrative example, for estimating the normal mean under $N(\mu, \sigma^2)$,  $Q=\bar X_n - \mu$ and $H=S(\mu|X)= n(\bar X_n - \mu)/\sigma^2=nG/\sigma^2$, and hence \eqref{eq:Rao} becomes equality,  whereas such an $H$ is clearly not permissible for \eqref{eq:GH}. 

Nevertheless, both inequalities reveal the tension between individual variations---features of their respective marginal distributions---and  their co-variation, which reflects their relationships, probabilistic or not. For \eqref{eq:Rao}, in order to keep $\Cov^2(G, H)$ at the value of $[Q^{\prime}(\theta)]^2>0$, the two variances $\Var(H)$ and $\Var(G)$ cannot be simultaneously small to an arbitrary degree, just as a rectangle cannot have arbitrarily small sides simultaneously when its area is bounded away from zero. This restriction leads to the Cram\'er-Rao lower bound. In \eqref{eq:Rao}, we purposefully write the Fisher information as the variance of the score function instead of the expectation of its negative derivative. The variance expression makes it clearer the co-variation essence of Cram\'er-Rao inequality, and draws a direct parallel with the inequality underlying HUP. 

Specifically, using the notation and the inequality \eqref{eq:hbound} of Section~\ref{sec:covoperator} and taking $G=\hat x$ and $H=\hat p$, we have 
\begin{equation}\label{eq:Heisenberg}
   \frac{\hbar^2}{4}\le |\Cov(G, H)|^2 \le \Var(H) \Var(G),
\end{equation}
 Comparing \eqref{eq:Heisenberg}
with \eqref{eq:Rao}, we see that the Cram\'er-Rao bound and the Heisenberg uncertainty principle are consequences of essentially the same statistical phenomena, that is, two marginal variances necessarily compete with each for being arbitrarily small, when the corresponding covariance is constrained in magnitude from below.

In contrast, for \eqref{eq:GH}, the trade-off is between the covariance and one of the marginal variances. To see this clearly, we can assume $\Var(H)=1$, which does not offend the assumption that $\E(H)=0$. Inequality \eqref{eq:GH} then becomes 
\begin{equation}\label{eq:GHs}
   \Cov^2(G, H) \le \Var(G) - R^{\rm opt} = R_G,
\end{equation}
where $R_G$ is the regret of $G$. On the surface, the changes of covariance and $\Var(G)$ appear to be coordinated instead of in competition, because the larger $\Cov^2(G, H)$, the larger $\Var(G)$. The reverse holds when the inequality is equality (which often is the case), and more broadly larger $\Var(G)$---and hence larger regret---at least allows more room for $\Cov^2(G, H)$ to grow.  But this is exactly where the tension lies when we want to improve both the learning and error assessment; improving learning means to reduce $R_G$ and hence have a \textit{smaller} $\Var(G)$, but improving error assessment requires a \textit{larger} $\Cov^2(G, H)$.  

\section{Elementary mathematics, advanced statistics, and inspiring philosophy}\label{sec:philosophy}

Mathematically, the proof of either \eqref{eq:Rao} or \eqref{eq:Heisenberg} is elementary, yet the implications of either inequality, as we know, are profound. Similarly, the inequality \eqref{eq:GH} is built upon equally elementary mathematics, and the work of \cite{bates2024cross} has already suggested its potential impact. However, many more studies remain, particularly regarding alternative loss functions, where the relevance of error assessment may not align with covariance. From a probabilistic standpoint, a thorough theoretical exploration of the relevance of an error assessor, $\hat{\delta}$, for the true error $\delta$ should involve investigating the joint distribution of $\hat{\delta}$ and $\delta$. In this context, irrelevance can be characterized by the independence between $\hat{\delta}$ and $\delta$.

On a broader level, formulating a general trade-off between learning and error assessment remains a complex task. This challenge stems from the need to define and measure the actual information utilized during learning and to identify relevant replications when assessing errors. Both `information' and `learning' are elusive notions, having taken on numerous interpretations throughout history, many of which require a refined understanding. For instance, even in the case of classical likelihood inference within parametric models, the role of conditioning in error assessment continues to provoke theoretical and practical debates.

I was reminded of this reality by an astrostatistics  project involving correcting conceptual and methodological errors in astrophysics for conducting model fitting and goodness-of-fit assessment via the popular C-statistics,  which is the likelihood ratio statistic under a Poisson regression model \citep{cash1979}. When the project started, I naively believed that it would be merely an exercise of applying classical likelihood theory and methods, perhaps with some clever computational tricks or approximations to render them practically efficient and hence appeals to astrophysicists.  

As reported in \cite{chen2024boosting}, however, the issue about whether one should condition on the MLE itself or not in the context of goodness-of-fit testing, is a rather nuanced one.  The issue is closely related to the issue of conditioning on ancillary statistics, since for testing distributional shape, the parametric parameters are \textit{nuisance objects} \citep[as termed in][]{meng2024bffer} and their MLE can be intuitively perceived as locally ancillary \citep{cox1980local,severini1993local} because the distribution shape of the MLE will be normal to the first order (under regularity conditions) despite the shape of the distribution being tested.  However, it is not exactly ancillary, and to decide when conditioning is beneficial (\textit{e.g.}, leading to a more powerful test) in any sample settling is not a straightforward matter. Higher order asymptotics can help provide insight, but communicating them intuitively is a tall order even for statisticians, let alone for astrophysicists or any scientists (including data scientists).  

However, regardless of whether low-level mathematics or high/tall order of statistics are involved, the ultimate challenge of contemplating and formulating uncertainty principles is epistemological, or even metaphysical.  For readers interested in philosophical contemplation---and I'd expect that statisticians should be in that group because statistics is essentially \textit{applied epistemology}\footnote{This was a characterization given by philosopher Hanti Lin during the JSM 2024, where Hanti and I co-organized a session where each philosopher presented for 20 minutes followed by a 15-min discussion by a statistician, and there were three pairs in total. (I made a mistake that embodied the statisticians' modesty: the estimated room size I provided to the JSM meeting department had an unacceptably negative bias.)}, I highly recommend the over 50 pages entry titled ``The Uncertainty Principle" by \cite{hilgevoord2024uncertainty} in \textit{The Stanford Encyclopedia of Philosophy}.\footnote{SEP is simply a fountain of afflatus and a \textit{Who's Who} in philosophy. Indeed SEP was where I came across Hanti Lin's  115-page entry on ``Bayesian Epistemology" \citep{lin_bayesian_2022}, and led to my invitation to Hanti to serve as a co-editor to establish the ``Meta Data Science" column ({\tt https://hdsr.mitpress.mit.edu/meta-data-science}) for \textit{Harvard Data Science Review}.} It is an erudite and thought-provoking  essay about the intellectual journey of Heisenberg's uncertainty principle. Even or perhaps especially the name ``uncertainty principle" has an interesting story behind it, because initially the name did not contain either `uncertainty' or `principle'.  

As \cite{hilgevoord2024uncertainty} discussed, the term \textit{uncertainty} has multiple meanings, and it is not obvious in which sense the phenomena revealed by \cite{heisenberg1927} qualifies as `uncertainty'; indeed, historically terms such as ``inaccuracy, spread,
imprecision, indefiniteness, indeterminateness, indeterminacy, latitude'' were used by various writers for what is now known as HUP.  More intriguingly, Heisenberg did not postulate the finding as any kind of \textit{principle}, but rather as \textit{relations}, such as ``inaccuracy relations" or ``indeterminacy relations". The discussions in Section~\ref{sec:covoperator} certainly reflect the relational nature of HUP, because it is fundamentally about the co-variation of position and momentum at the mechanism level. 

The entry by \cite{hilgevoord2024uncertainty} invites readers to consider a fundamental question that underpins these onomasiological reflections: Is the HUP a mere epistemic constraint, or a metaphysical limitation in nature? Unsurprisingly, this question is a source of ongoing dispute among philosophers of physics and even among physicists themselves.
The most well-known historical debates are Heisenberg and Bohr's Copenhagen interpretation emphasizing the metaphysical indeterminacy, and the contrasting deterministic interpretation developed by de Broglie and Bohm, known as Bohmian mechanics \citep{hilgevoord2024uncertainty}.

Given I have already greatly exceeded the deadline to submit this  essay,  I will refrain from revealing any further thrills provided in \cite{hilgevoord2024uncertainty}, such as more recent debates about HUP,  leaving readers to enjoy their own treasure hunt.  But I will mention that this question has prompted me to wonder whether inequality \eqref{eq:corrinq} also suggests that any effort to assess the actual error is antithetic to probabilistic learning. 

This is because the crux of probabilistic learning—unlike deterministic approaches, such as solving algebraic equations—lies in using distributions as our fundamental mathematical vehicles for carrying our states of knowledge (or lack thereof) and for transporting data into information that furthers learning. From this distributional perspective, assessing the actual error means to assess the \textit{distribution} of the actual error, which is all we need to, for example,  provide the usual confidence regions. It does suffer from the leap of faith problem as discussed in Section~\ref{sec:Jay}, but then that is a universal predicament to any form of empirical learning, as far as I can imagine. 
 
\section{From uncertainty principles to happy marriages ...}\label{sec:marriage}

A further inspiration from \cite{hilgevoord2024uncertainty} is its discussion on the relationship between the original semi-quantitative argument made by \cite{heisenberg1927} and the mathematical formalism established by \cite{kennard1927}. Kennard's inequality \eqref{eq:keen} is precise, but can be perceived being narrow, for instance, in its reliance on standard deviation to describe ``uncertainty." A similar limitation applies to inequality \eqref{eq:corrinq}, which assesses relevance through linear correlation, a measure surely is not universally appropriate for capturing the notion of relevance. 

More broadly, much remains to be examined regarding the trade-offs between the flexibility of qualitative frameworks, which embrace the nuances and ambiguities of natural language, and the rigor of quantitative formulations, which offer the precision of mathematical language but often at the risk of being overly restrictive or idealized. Reflecting on these trade-offs is essential to learning. Statisticians and data scientists, in particular, can draw from centuries of philosophical inquiry into epistemology, as exemplified by the discussions surrounding the HUP and the like. 

In truth, when thoughtfully practiced, data science embodies—or ought to embody—a harmonious blend of quantitative and qualitative thinking and reasoning. This was the central theme of my \textit{Harvard Data Science Review} editorial, “Data Science: A Happy Marriage of Quantitative and Qualitative Thinking?” \citep{Meng2021Data}, inspired by \cite{Tanweer2021Why}’s compelling article, “Why the Data Revolution Needs Qualitative Thinking.” Maintaining this harmony, akin to sustaining a functioning marriage, requires commitment from all parties and a willingness to compromise. Ultimately, it calls for the wisdom to recognize that individual fulfillment and happiness—whether in marriage, mentorship, or mind melding or mating—depends profoundly on collective well-being. 
Professor Rao certainly embodied this wisdom.

I vividly recall my first visit to Pennsylvania State University as a seminar speaker, shortly after Professor Rao’s 72nd birthday on September 10, 1992. During the seminar lunch, Professor Rao graciously joined us. We—students and early-career researchers (myself included, back when my hair was dense almost surely everywhere)—felt honored by his presence. All questions naturally revolved around statistics, except for one that made us all chuckle: “Professor Rao, how does one live a long and happy life?”

Without missing a beat, and with his characteristic paced, confident cadence, Rao replied, “Keep your wife happy.”

\section{A prologue or an invitation}\label{sec:prologue}

For those who would like this article to conclude with a statistical Q\&A: During the elevator ride following my seminar, which carried the seemingly oxymoronic title “A Bayesian p-value” (a deliberate contrast to the title of \cite{meng1994posterior}), Professor Rao turned to me and asked, “Do people still use \textit{p}-values?” To which I responded…

Well, I’ll leave that as a missing data point, inviting you to impute your own favorite answer. Alternatively, if you prefer, find a deliberately embedded mathematical (but petty) error in this article and exchange it for the answer by emailing {\tt meng@stat.harvard.edu} (as long as God permits me to respond). 

\section*{Acknowledgments}

I am deeply grateful to physicists Aurore Courtoy, Louis Lyons, Thomas Junk, and Pavel Nadolsky, as well as statistician Yazhen Wang, for their careful and patient explanations regarding the non-existence of a joint probabilistic distribution of a particle's position and momentum. I am equally indebted to Hanti Lin for elucidating the philosophical debates surrounding the Heisenberg Uncertainty Principle.

My thanks also extend to editor Bhramar Mukherjee, to whom I owe a profound debt, and to Peter Bickel, Joe Blitzstein, Radu Craiu, Walter Dempsey, Benedikt Höltgen, Peter McCullagh, Pavlos Msaouel, Steve Stigler, Robert Tibshirani, Théo Voldoire, and Bob Williamson, for collectively providing insightful comments and sharing relevant literature—some of which may inspire a sequel to this essay.

I also thank Julie Vu and Sicheng Zhou for their meticulous proofreading efforts; naturally, any remaining errors are entirely my own (though I wish they weren’t!). Finally, I acknowledge partial financial support from the NSF during the period when this essay was conceived and completed.

\section*{Appendix A:  Derivations for The Regression Example in Section~\ref{sec:lunch}}
In general, the weighted estimate of $\theta$ can be written as 
\begin{equation*}
\hat\theta_w = \frac{\sum_{i=1}^n w_i X_iY_i}{\sum_{i=1}^n w_i X_i^2} ,
\end{equation*}
with OLS corresponding to choosing $w_i=1$ and BLUE given by $w_i=\sigma^{-2}_i$, for all $i$. Conditioning on $\mathbf X$ but for notational simplicity we suppress the conditioning notation in all expectations below, we have 
\begin{equation*}
\Var(\hat\theta_w) =\frac{\sum_{i=1}^nw_i^2X_i^2\sigma_i^2}{[\sum_{i=1}^nw_i X_i^2]^2}=\frac{T_{w,\sigma}}{T_w^2}.
\end{equation*}
Let $\hat{r}_{w,j}=Y_j-\hat\theta_w X_j$. Because $\E(\hat{r}_{w,j})=0$, to calculate $\rho$, we only need to calculate
\begin{align*}\label{eq:OSL}
\E[\hat\theta_w(Y_j-\hat\theta_wX_j)] &= \frac{\sum_{i=1}^n w_iX_i\E[Y_iY_j]}{T_w} - \frac{\E[\sum_{i=1}^n w_iX_iY_i]^2 X_j}{T^2_w}\\
&= \frac{\sum_{i=1}^n w_iX_i[\Cov(Y_i, Y_j)+\theta^2 X_iX_j]}{T_w} - \frac{\left[\sum_{i=1}^n w_i^2X_i^2\sigma^2_i + \theta^2 T^2_w \right]X_j}{T^2_w}\\
&= \frac{(\theta^2 T_w +w_j\sigma^2_j)X_j}{T_w} - \frac{\left[T_{w,\sigma} + \theta^2 T^2_w \right]X_j}{T^2_w} = \frac{X_j}{T_w} \left[w_j\sigma^2_j - \frac{T_{w,\sigma}}{T_w} \right];
\end{align*}
and 
\begin{align*}
\Var(\hat{r}_{w,j})&= \Var\left[\frac{\sum_{i=1}^n w_iX_i (X_iY_j - X_j Y_i)}{T_w}\right]=T^{-2}_w \Var\left[\sum_{i\not=j}^n w_iX_i (X_iY_j - X_j Y_i)\right] \\
&=T^{-2}_w \E\left\{\Var\left[\sum_{i\not=j}^n w_iX_i (X_iY_j - X_j Y_i)|Y_j\right]\right\}+ \Var\left\{\E\left[\sum_{i\not=j}^n w_iX_i (X_iY_j - X_j Y_i)|Y_j\right]\right\}\\
&=T^{-2}_w\left\{\left[X_j^2\sum_{i\not=j}^n w_i^2X_i^2\sigma^2_i\right]+ \Var\left[\sum_{i\not=j}^nw_i X_i^2 Y_j\right]\right\}\\
&=T^{-2} _w\left\{\left[X_j^2(T_{w,\sigma} -w_j^2 X_j^2\sigma^2_j)\right]+ [T_w - w_jX_j^2]^2 \sigma^2_j\right\}\\
&=T^{-2} _w\left\{X_j^2T_{w,\sigma} + \sigma^2_j[T_w^2 - 2T_ww_jX^2_j]\right\}.
\end{align*}
Putting all the pieces together, we have 
\begin{equation}\label{eq:corrwls}
\Corr(\hat\theta_w, \hat{r}_{w,j}) 
= \frac{X_j\left(w_j\sigma^2_jT_w - T_{w,\sigma}\right)}{
{\sqrt{T_{w, \sigma}\left[X_j^2T_{w,\sigma} + \sigma^2_j(T_w^2 - 2T_ww_jX^2_j)\right] }}}, \quad j=1, 2.
\end{equation}
For $n=2, j=1$, expression \eqref{eq:corrwls} simplifies to the desired \eqref{eq:corrn2} because 
\begin{align*}
\Corr(\hat\theta_w, r_{w,1}) 
&= \frac{X_1X_2^2w_2(w_1\sigma^2_1-w_2\sigma^2_2)}{\sqrt{[X_1^2w^2_2X_2^2\sigma^2_2+ w_2^2X_2^4\sigma^2_1][w_1^2X_1^2 \sigma^2_1+w^2_2X_2^2 \sigma^2_2]}}\\
&=\frac{X_1|X_2|(w_1\frac{\sigma_1}{\sigma_2}-w_2\frac{\sigma_2}{\sigma_1})}{\sqrt{[X_1^2\sigma^{-2}_1+ X_2^2\sigma^{-2}_2][w_1^2X_1^2 \sigma^2_1+w^2_2X_2^2 \sigma^2_2]}}.
\end{align*}
To calculate the relative regret (RR), we have 
\begin{equation}\label{eq:corra}
\Var(\hat\theta_w)=\Var\left[\frac{\sum_{i=1}^nw_i X_i Y_j}{T_w}\right] =\frac{w_1^2X_1^2 \sigma^2_1+w_2^2X_2^2 \sigma^2_2}{[w_1X_1^2+w_2X_2^2]^2},
\end{equation}
which also implies, by taking $w_i\propto \sigma^{-2}_i$,
\begin{equation}\label{eq:bluea}
\Var(\hat\theta_{\rm BLUE})=\frac{1}{(X_1^2 \sigma^{-2}_1+X_2^2 \sigma^{-2}_2)}.
\end{equation}
Putting together \eqref{eq:corra} and \eqref{eq:bluea} yields the desired \eqref{eq:regret}.

\section*{Appendix B: Derivation of \eqref{eq:joincorr} in Section~\ref{sec:Jay}}

Because $\hat\delta^2$ and $\delta^2$ are independent given $\theta=\{\mu, \sigma^2\}$ and hence $\Cov(\hat\delta^2, \delta^2|\mu, \sigma^2)=0$, we see over the joint replication, 
\begin{align*}
    \Cov(\hat\delta^2, \delta^2)=\E\left[\Cov(\hat\delta^2, \delta^2|\mu, \sigma^2)\right] + \Cov\left[\E(\hat\delta^2|\mu, \sigma^2), \E(\delta^2|\mu, \sigma^2)\right]=\frac{1}{n^2}\Var(\sigma^2),
    \end{align*}
    as long as the prior distribution for $\theta=\{\mu, \sigma^2\}$ is proper. Furthermore, conditioning on $\theta=\{\mu, \sigma^2\}$,  $\delta^2\sim \sigma^2\chi^2_1/n$ and $\hat\delta^2\sim \sigma^2\chi^2_{n-1}/[n(n-1)]$ (where the two chi-square variables are independent of each other), we have 
    \begin{align*}\label{eq:jointvar}
     \Var(\hat\delta^2)=&\E\left[\Var(\hat\delta^2|\mu, \sigma^2)\right] + \Var\left[\E(\hat\delta^2|\mu, \sigma^2)\right]= \frac{2}{(n-1)n^2}\E\left(\sigma^4\right)+\frac{1}{n^2}\Var(\sigma^2);\\
     \Var(\delta^2)=&\E\left[\Var(\delta^2|\mu, \sigma^2)\right] + \Var\left[\E(\delta^2|\mu, \sigma^2)\right]= \frac{2}{n^2}\E\left(\sigma^4\right)+\frac{1}{n^2}\Var(\sigma^2).
\end{align*}
Consequently, we see over the joint replication, 
\begin{align*}
    \Corr(\hat\delta^2, \delta^2)=\frac{\Var(\sigma^2)}{\sqrt{2(n-1)^{-1}\E(\sigma^4)+\Var(\sigma^2)}\sqrt{2\E(\sigma^4)+\Var(\sigma^2)}},
\end{align*}
which yields \eqref{eq:joincorr} because $\E(\sigma^4)=\Var(\sigma^2)+[\E(\sigma^2)]^2$.

\section*{Appendix C: A quasi-score analogy for understanding the lack of joint probability}\label{apx:score}

For statistically oriented readers, an instructive—though far from being perfect—analogy to the issue of the non-existence of a probabilistic model due to violations of symmetry or commutativity is the generalization from likelihood inference via the score function to estimation based on quasi-score functions. The correct score function, when available, provides the most efficient inference asymptotically (under regularity conditions). However, specifying the correct data-generating model often requires more information and resources than we typically possess.

In contrast, a quasi-score function only requires the specification of the first two moments of the data-generating model. This makes it a more practical and robust alternative to exact model-based inference, particularly in the presence of model misspecification. However, this robustness comes at the cost of reduced efficiency, reflecting the trade-off inherent in this approach.

Broadly speaking there are three types of pseudo scores: (I) those that are equivalent to the actual score; (II) those that are not equivalent to the actual score, but are equivalent to the score from a misspecified data generating model, and (III) those that cannot be derived from any probabilistic model.  

Type (III) exists because any (differentiable) authentic score vector $\left(S_1(\theta), \ldots, S_d(\theta)\right)^\top$ for a $d$-dimension parameter $\theta=\left(\theta_1, \ldots, \theta_d\right)^\top$ must satisfy 
\begin{equation}\label{eq:scorecond}
\frac{\partial S_i(\theta)}{\partial \theta_j}= \frac{\partial S_j(\theta)}{\partial \theta_i}, \quad \forall\ i, j =1, \ldots, d,
\end{equation}
because the corresponding (observed) Fisher information matrix, $-\frac{\partial S(\theta)}{\partial \theta}$, is symmetric. However, even some most innocent looking quasi-scores, such as for certain $2\times 2$ contingency tables, the symmetry requirement of \eqref{eq:scorecond} can be easily violated, as demonstrated in Chapter 9 of \cite{mccullagh1989generalized}, which is an excellent source for understanding quasi scores and estimation equations in general. 

The fact that violating the symmetry condition \eqref{eq:scorecond} rules out the possibility of being an actual score may help some of us imagine how the lack of symmetry or commutativity might rule out the existence of a probability specification, at least from a mathematical perspective. Furthermore, just as one can generalize from likelihood to quasi-likelihood of many shapes and forms---again see \cite{mccullagh1989generalized}---the non-existence of a probabilistic distribution does not prevent us from forming quasi-distributions for various purposes, such as the Wigner quasiprobability distribution, which permits negative values, for position and momentum $(x, p)$ \citep{hillery1984distribution,lorce2011quark}. Whether the mechanism-level covariances as given in \eqref{eq:covcommon}-\eqref{eq:covcommonc} have the same magnitude as that from the Wigner quasiprobability distribution will be left as a homework exercise. 
\bibliographystyle{sa}
\bibliography{ref}

\begin{thebibliography}{}

\bibitem[Abba\textit{ et~al.}, 2024]{abba2024bayesian}
Abba, M.~A., Williams, J.~P., and Reich, B.~J. (2024).
\newblock A {Bayesian} shrinkage estimator for transfer learning.
\newblock {\em arXiv:2403.17321}, {\bfseries }.

\bibitem[Bates\textit{ et~al.}, 2024]{bates2024cross}
Bates, S., Hastie, T., and Tibshirani, R. (2024).
\newblock Cross-validation: what does it estimate and how well does it do it?
\newblock {\em Journal of the American Statistical Association}, {\bfseries 119}, 1434--1445.

\bibitem[Berger\textit{ et~al.}, 2024]{berger2024handbook}
Berger, J., Meng, X.-L., Reid, N., and Xie, M.-g. (2024).
\newblock {\em Handbook of Bayesian, Fiducial, and Frequentist Inference}.
\newblock CRC Press.

\bibitem[Blitzstein and Hwang, 2014]{blitzstein_hwang_2014}
Blitzstein, J.~K. and Hwang, J. (2014).
\newblock {\em Introduction to Probability}.
\newblock CRC Press, Boca Raton, FL, 1st edition.

\bibitem[Casella and Berger, 2024]{casella2024statistical}
Casella, G. and Berger, R. (2024).
\newblock {\em {Statistical Inference}}.
\newblock CRC Press.

\bibitem[{Cash}, 1979]{cash1979}
{Cash}, W. (1979).
\newblock Parameter estimation in astronomy through application of the likelihood ratio.
\newblock {\em The Astrophysical Journal}, {\bfseries 228}, 939.

\bibitem[Chapman and Robbins, 1951]{chapman1951minimum}
Chapman, D.~G. and Robbins, H. (1951).
\newblock Minimum variance estimation without regularity assumptions.
\newblock {\em The Annals of Mathematical Statistics}, {\bfseries 22}, 581--586.

\bibitem[Chen\textit{ et~al.}, 2024]{chen2024boosting}
Chen, Y., Li, X., Meng, X.-L., van Dyk, D.~A., Bonamente, M., and Kashyap, V. (2024).
\newblock Boosting {C-statistics} in astronomy via conditioning: More power, less computation.
\newblock Technical report, Department of Statistics, University of Michigan.

\bibitem[Cox, 1980]{cox1980local}
Cox, D.~R. (1980).
\newblock Local ancillarity.
\newblock {\em Biometrika}, {\bfseries 67}, 279--286.

\bibitem[Craiu\textit{ et~al.}, 2023]{craiu2023six}
Craiu, R.~V., Gong, R., and Meng, X.-L. (2023).
\newblock Six statistical senses.
\newblock {\em Annual Review of Statistics and Its Application}, {\bfseries 10}, 699--725.

\bibitem[Dembo, 1990]{dembo1990information}
Dembo, A. (1990).
\newblock Information inequalities and uncertainty principles.
\newblock {\em Department of Statistics, Stanford University., Stanford, CA, Technical Report}, {\bfseries 75}.

\bibitem[Dembo\textit{ et~al.}, 1991]{dembo1991information}
Dembo, A., Cover, T.~M., and Thomas, J.~A. (1991).
\newblock Information inequalities and uncertainty principles.
\newblock {\em IEEE Transactions on Information Theory}, {\bfseries 37}, 1501--1518.

\bibitem[Dempster, 1963]{Dempster1963c}
Dempster, A.~P. (1963).
\newblock Further examples of inconsistencies in the fiducial argument.
\newblock {\em The Annals of Mathematical Statistics}, {\bfseries 34}, 884--891.

\bibitem[Efron, 2020]{efron2020prediction}
Efron, B. (2020).
\newblock Prediction, estimation, and attribution.
\newblock {\em International Statistical Review}, {\bfseries 88}, S28--S59.

\bibitem[Efron and Hinkley, 1978]{efron1978assessing}
Efron, B. and Hinkley, D.~V. (1978).
\newblock Assessing the accuracy of the maximum likelihood estimator: Observed versus expected fisher information.
\newblock {\em Biometrika}, {\bfseries 65}, 457--483.

\bibitem[Gelman and Betancourt, 2013]{gelman2013does}
Gelman, A. and Betancourt, M. (2013).
\newblock Does quantum uncertainty have a place in everyday applied statistics.
\newblock {\em Behavioral and Brain Sciences}, {\bfseries 36}, 285.

\bibitem[Gong and Meng, 2021]{gong2021judicious}
Gong, R. and Meng, X.-L. (2021).
\newblock Judicious judgment meets unsettling updating: Dilation, sure loss, and simpson’s paradox.
\newblock {\em Statistical Science}, {\bfseries 36}, 169--214.
\newblock Discussion article with rejoinder.

\bibitem[Griffiths and Schroeter, 2018]{griffiths2018introduction}
Griffiths, D.~J. and Schroeter, D.~F. (2018).
\newblock {\em {Introduction to Quantum Mechanics}}.
\newblock Cambridge University Press.

\bibitem[Heisenberg, 1927]{heisenberg1927}
Heisenberg, W. (1927).
\newblock Über den anschaulichen inhalt der quantentheoretischen kinematik und mechanik.
\newblock {\em Zeitschrift für Physik}, {\bfseries 43}, 172--198.

\bibitem[Hilgevoord and Uffink, 2024]{hilgevoord2024uncertainty}
Hilgevoord, J. and Uffink, J. (2024).
\newblock The uncertainty principle.
\newblock In Zalta, E.~N. and Nodelman, U., editors, {\em The Stanford Encyclopedia of Philosophy}. Stanford University.
\newblock Spring 2024 edition.

\bibitem[Hillery\textit{ et~al.}, 1984]{hillery1984distribution}
Hillery, M., O'Connell, R.~F., Scully, M.~O., and Wigner, E.~P. (1984).
\newblock Distribution functions in physics: Fundamentals.
\newblock {\em Physics reports}, {\bfseries 106}, 121--167.

\bibitem[Hoeffding, 1940]{Hoeffding1940}
Hoeffding, W. (1940).
\newblock Ma{\ss}tabinvariante {K}orrelatiostheorie.
\newblock {\em Schriften des Mathematischen Instituts und des Instituts f\"ur Angewandte Mathematik der Universit\"{a}t Berlin}, {\bfseries 5}, 179--233.

\bibitem[Kennard, 1927]{kennard1927}
Kennard, E.~H. (1927).
\newblock Zur quantenmechanik einfacher bewegungstypen.
\newblock {\em Zeitschrift für Physik}, {\bfseries 44}, 326--352.

\bibitem[Landau and Lifshitz, 2013]{landau2013quantum}
Landau, L.~D. and Lifshitz, E.~M. (2013).
\newblock {\em {Quantum Mechanics: Non-relativistic Theory}}, volume~3.
\newblock Elsevier.

\bibitem[Le~Cam, 1956]{le1956asymptotic}
Le~Cam, L. (1956).
\newblock On the asymptotic theory of estimation and testing hypotheses.
\newblock In {\em Proceedings of the Third Berkeley Symposium on Mathematical Statistics and Probability, Volume 1: Contributions to the Theory of Statistics}, volume~3, pages 129--157. University of California Press.

\bibitem[Lehmann and Casella, 2006]{lehmann2006theory}
Lehmann, E.~L. and Casella, G. (2006).
\newblock {\em {Theory of Point Estimation}}.
\newblock Springer Science \& Business Media.

\bibitem[Lin, 2024a]{lin_bayesian_2022}
Lin, H. (2024a).
\newblock Bayesian epistemology.
\newblock In Zalta, E.~N. and Nodelman, U., editors, {\em The Stanford Encyclopedia of Philosophy}.
\newblock 2024 Edition, originally published 2022.

\bibitem[Lin, 2024b]{Lin2024To}
Lin, H. (2024b).
\newblock To be a {Frequentist} or {Bayesian}? {Five} {positions} in a {spectrum}.
\newblock {\em Harvard Data Science Review}, {\bfseries 6}.
\newblock https://hdsr.mitpress.mit.edu/pub/axvcupj4.

\bibitem[Liu and Meng, 2014]{liu2014comment}
Liu, K. and Meng, X.-L. (2014).
\newblock Comment: A fruitful resolution to simpson’s paradox via multiresolution inference.
\newblock {\em The American Statistician}, {\bfseries 68}, 17--29.

\bibitem[Liu and Meng, 2016]{liu2016there}
Liu, K. and Meng, X.-L. (2016).
\newblock There is individualized treatment. {Why} not individualized inference?
\newblock {\em Annual Review of Statistics and Its Application}, {\bfseries 3}, 79--111.

\bibitem[Lorce and Pasquini, 2011]{lorce2011quark}
Lorce, C. and Pasquini, B. (2011).
\newblock Quark wigner distributions and orbital angular momentum.
\newblock {\em Physical Review D—Particles, Fields, Gravitation, and Cosmology}, {\bfseries 84}, 014015.

\bibitem[McCullagh, 1999]{mccullagh1999discussion}
McCullagh, P. (1999).
\newblock Discussion on {Lindsey, J.K. (1999). "Some statistical heresies"}.
\newblock {\em Journal of the Royal Statistical Society: Series D (The Statistician)}, {\bfseries 48}, 34--35.

\bibitem[McCullagh and Nelder, 1989]{mccullagh1989generalized}
McCullagh, P. and Nelder, J.~A. (1989).
\newblock {\em Generalized Linear Models}, volume~37 of {\em Monographs on Statistics and Applied Probability}.
\newblock Chapman \& Hall/CRC, London, 2nd edition.

\bibitem[Meng, 1994]{meng1994posterior}
Meng, X.-L. (1994).
\newblock Posterior predictive $ p $-values.
\newblock {\em The annals of statistics}, {\bfseries 22}, 1142--1160.

\bibitem[Meng, 2018]{meng2018statistical}
Meng, X.-L. (2018).
\newblock Statistical paradises and paradoxes in big data {(I)}: law of large populations, big data paradox, and the 2016 us presidential election.
\newblock {\em The Annals of Applied Statistics}, {\bfseries 12}, 685--726.

\bibitem[Meng, 2021]{Meng2021Data}
Meng, X.-L. (2021).
\newblock Data science: A happy marriage of quantitative and qualitative thinking?
\newblock {\em Harvard Data Science Review}, {\bfseries 3}.
\newblock https://hdsr.mitpress.mit.edu/pub/pger71uh.

\bibitem[Meng, 2024]{meng2024bffer}
Meng, X.-L. (2024).
\newblock {A BFFer’s exploration with nuisance constructs: Bayesian \textit{p}-value, H-likelihood, and Cauchyanity.}
\newblock In {\em Handbook of Bayesian, Fiducial, and Frequentist Inference, Eds J. Berger, XL. Meng, N. Reid and M. Xie}, pages 161--187. Chapman and Hall/CRC.

\bibitem[Rao, 1945]{rao1945information}
Rao, C.~R. (1945).
\newblock Information and the accuracy attainable in the estimation of statistical parameters.
\newblock {\em Bulletin of the Calcutta Mathematical Society}, {\bfseries 37}, 81--91.

\bibitem[Rao, 1962]{rao1962efficient}
Rao, C.~R. (1962).
\newblock Efficient estimates and optimum inference procedures in large samples.
\newblock {\em Journal of the Royal Statistical Society: Series B (Methodological)}, {\bfseries 24}, 46--63.

\bibitem[Severini, 1993]{severini1993local}
Severini, T.~A. (1993).
\newblock Local ancillarity in the presence of a nuisance parameter.
\newblock {\em Biometrika}, {\bfseries 80}, 305--320.

\bibitem[Stam, 1959]{stam1959some}
Stam, A.~J. (1959).
\newblock Some inequalities satisfied by the quantities of information of fisher and shannon.
\newblock {\em Information and Control}, {\bfseries 2}, 101--112.

\bibitem[Tanweer\textit{ et~al.}, 2021]{Tanweer2021Why}
Tanweer, A., Gade, E.~K., Krafft, P., and Dreier, S. (2021).
\newblock Why the data revolution needs qualitative thinking.
\newblock {\em Harvard Data Science Review}, {\bfseries 3}.
\newblock https://hdsr.mitpress.mit.edu/pub/u9s6f22y.

\bibitem[T{\'o}th and Fr{\"o}wis, 2022]{toth2022uncertainty}
T{\'o}th, G. and Fr{\"o}wis, F. (2022).
\newblock Uncertainty relations with the variance and the quantum fisher information based on convex decompositions of density matrices.
\newblock {\em Physical Review Research}, {\bfseries 4}, 013075.

\bibitem[T{\'o}th and Petz, 2013]{toth2013extremal}
T{\'o}th, G. and Petz, D. (2013).
\newblock Extremal properties of the variance and the quantum fisher information.
\newblock {\em Physical Review A—Atomic, Molecular, and Optical Physics}, {\bfseries 87}, 032324.

\bibitem[Wald, 1943]{wald1943tests}
Wald, A. (1943).
\newblock Tests of statistical hypotheses concerning several parameters when the number of observations is large.
\newblock {\em Transactions of the American Mathematical society}, {\bfseries 54}, 426--482.

\bibitem[Wang, 2022]{Wang2022When}
Wang, Y. (2022).
\newblock When quantum computation meets data science: {Making} data science quantum.
\newblock {\em Harvard Data Science Review}, {\bfseries 4}.
\newblock https://hdsr.mitpress.mit.edu/pub/kpn45eyx.

\end{thebibliography}

\end{document}